\documentclass[reqno,a4paper,12pt]{amsart} 

\usepackage{amsmath,amscd,amsfonts,amssymb}
\usepackage{mathrsfs,dsfont}

\usepackage{tikz}
\usetikzlibrary{patterns}
\usepackage{float}
\usepackage[toc]{appendix}
\usepackage{caption}
\usepackage{subcaption}
\usepackage{enumerate}

\numberwithin{equation}{section}
\numberwithin{figure}{section}

\def\R{\mathbb{R}}
\def\C{\mathbb{C}}

\def\Z{\mathbb{Z}}

\def\M{\mathcal{M}}

\def\1{\mathds{1}}

\def\dH{\dim_{\mathcal{H}}}

\renewcommand\leq{\leqslant}
\renewcommand\geq{\geqslant}

\renewcommand\hat{\widehat}
\renewcommand\Re{\operatorname{Re}}
\renewcommand\Im{\operatorname{Im}}

\newcommand{\supp}{\operatorname{supp}}
\newcommand{\sgn}{\operatorname{sgn}}
\newcommand{\Span}{\operatorname{Span}}

\theoremstyle{plain}
\newtheorem{thm}{Theorem}[section]

\newtheorem{corollary}[thm]{Corollary}
\newtheorem{prop}[thm]{Proposition}

\newtheorem*{claim*}{Claim}
\newtheorem*{thm*}{Theorem}

\theoremstyle{definition}
\newtheorem{definition}[thm]{Definition}
\newtheorem*{definition*}{Definition}
\newtheorem*{remarks*}{Remarks}
\newtheorem*{remark*}{Remark}

\newenvironment{enumerate-math}
{\begin{enumerate}
\addtolength{\itemsep}{5pt}
}
{\end{enumerate}}

\newenvironment{enumerate-text}
{\begin{enumerate}
\addtolength{\itemsep}{5pt}
}
{\end{enumerate}}

\begin{document}

\title[Mixed-norm of orthogonal projections and analytic interpolation]{Mixed-norm of orthogonal projections and analytic interpolation on dimensions of measures}

\author{Bochen Liu}
\address{Department of Mathematics \& International Center for Mathematics, Southern University of Science and Technology, Shenzhen 518055, China}

\email{Bochen.Liu1989@gmail.com}

\date{}


\begin{abstract}
Suppose $\mu, \nu$ are compactly supported Radon measures on $\mathbb{R}^d$ and $V\in G(d,n)$ is an $n$-dimensional subspace. In this paper we systematically study the mixed-norm
$$\int\|\pi^y\mu\|_{L^p(G(d,n))}^q\,d\nu(y),\ \forall\,p,q\in[1,\infty),$$
where $\pi_V:\mathbb{R}^d\rightarrow V$ denotes the orthogonal projection and
$$\pi^y\mu(V)=\int_{y+V^\perp}\mu\,d\mathcal{H}^{d-n}=\pi_V\mu(\pi_Vy),\ \text{if $\mu$ has continuous density}.$$
When $n=d-1$ and $p=q$, our result significantly improves a previous result of Orponen.

In the proof we consider integer exponents first, then interpolate analytically, not only on $p,q$, but also on dimensions of measures. We also introduce a new quantity called $s$-amplitude, to present our results and illustrate our ideas. This mechanism provides new perspectives on operators with measures, thus has its own interest. We also give an alternative proof of a recent result of Dabrowski, Orponen, Villa on $\|\pi_V\mu\|_{L^p(\mathcal{H}^n\times G(d,n))}$.

The following consequences are also interesting.
\begin{itemize}
	\item We discover jump discontinuities in the range of $p$ at the critical line segment
	$$\{(s_\mu, s_\nu)\in(0,d)^2: s_\mu+s_\nu=2n,\, 0<s_\nu<n\},$$ 
	where $s_\mu, s_\nu$ are Frostman exponents of $\mu, \nu$ respectiely. This is unexpected and surprising.

	\item Given $1\leq m\leq d-1$ and $E, F\subset\R^d$, we obtain dimensional threshold on  whether there exists $y\in F$ such that
	$$\gamma_{d,m}\{V\in G(d,m): V=\Span\{x_1-y,\dots,x_m-y\}: x_1,\dots,x_m\in E\}>0.$$
	This generalizes the visibility problem ($m=1$). In particular,  when $m>\frac{d}{2}$ and $\dH E$ is large enough, the exceptional set has Hausdorff dimension $0$. 
\end{itemize}
\end{abstract}

\maketitle
\setcounter{tocdepth}{1}
\tableofcontents

\section{Introduction}
Let $\pi_e(x)=x\cdot e$, $x\in\R^2$, $e\in S^1$, denote the orthogonal projection. In 1954 Marstrand \cite{Mar54} proved that, for every $E\subset\R^2$, $\dH E>1$, the set $\pi_e(E)$ has positive Lebesgue measure for almost all $e\in S^1$. In 1968, Kaufman \cite{Kau68} gave a simple alternative proof of Marstrand projection theorem using Fourier analysis. Moreover, he proved that the induced measure $\pi_e\mu$ on $\pi_e(E)$ has $L^2$ density for almost all $e\in S^1$, where $\mu$ is a Frostman measure on $E$. Nowadays orthogonal projection has become a central problem in geometric measure theory, and has been studied actively from different perspectives. 

In higher dimensions, we denote orthogonal projections by $\pi_V:\R^d\rightarrow V$, where $V\in G(d,n)$ is a $n$-dimensional subspace of $\R^d$ and $G(d,n)$ denotes the Grassmannian. Also throughout this paper $\mathcal{H}^s$ denotes the $s$-dimensional Hausdorff measure, and $\gamma_{d,n}$ denotes the probability measure on $G(d,n)$ induced by the Haar measure $\theta_d$ on the orthogonal group.

Of course $\pi_V\mu\in L^2(\mathcal{H}^n|_V)$ is not a necessary condition for $\mathcal{H}^n(\pi_V (E))>0$. But for technical reasons the $L^2$-method is the most popular approach to this type of problems. We refer to \cite{Mat87}\cite{Liu18} for $L^2$-approaches to Falconer distance conjecture.

It is not surprising that $L^2$ does not always help. Let
$$\pi^y(x):=\frac{x-y}{|x-y|}\in S^{d-1}$$
denote the radial projection. The visibility problem asks that, given $E, F\subset\R^d$, whether there exists $y\in F$ such that $\pi^y(E)\subset S^{d-1}$ has positive surface measure. Although there are $L^2$-estimates in the literature (e.g., \cite{PS00}), their geometric consequences are far from being desired. Finally, Orponen \cite{Orp19} proved that, if $s>d-1$ and $s+t>2(d-1)$, then
\begin{equation}
	\label{Lp-Orponen}
	\int \|\pi^y\mu\|_{L^p(S^{d-1})}^p\,d\nu(y)\lesssim I_{s}(\mu)^{p/2}\cdot I_{t}(\nu)^{1/2}
\end{equation}
for every $$1\leq p<\min\left\{\frac{t}{2(d-1)-s},\ 2-\frac{t}{d-1}\right\},$$
where $I_s$ denotes the $s$-energy and $\pi^y\mu$ denotes the induced measure on $\pi^y(E)$, which equals $$\int_{y+e^\perp}\mu\,d\mathcal{H}^{d-1}=\pi_e\mu(\pi_ey)$$
if $\mu$ has continuous density.

As a geometric consequence of \eqref{Lp-Orponen}, if $E, F\subset\R^d$, $\dH E>d-1$, $\dH E+\dH F>2(d-1)$, then
$$\mathcal{H}^{d-1}(\pi^y (E))>0, \ \text{for some }y\in F.$$
It is known that neither of assumptions $\dH E>d-1$, $\dH E+\dH F>2(d-1)$ can be relaxed. We refer to \cite{Orp18} for the discussion on the sharpness and \cite{Mat15}, Example 5.13, for details of examples. This also implies both $s>(d-1)$ and $s+t>2(d-1)$ are necessary for the existence of $p>1$ in \eqref{Lp-Orponen}. 

What's more, the estimate \eqref{Lp-Orponen} was quickly introduced into the study of Falconer distance conjecture, playing crucial roles in a couple of recent breakthroughs \cite{KS18}\cite{GIOW18}. This brings more attention to $L^p$-estimates of projections.

Even more recently, Dabrowski, Orponen, Villa \cite{DOV22} proved that, if $\mu$ is a compactly supported measure on $\R^d$ satisfying the $s$-dimensional Frostman condition
$$\mu(B(x,r))\lesssim r^s,\ \forall\,r>0,\,x\in\R^d,$$
then
\begin{equation}
	\label{DOV}
	\int\|\pi_V\mu\|_{L^p(\mathcal{H}^n)}^p\,d\gamma_{d,n}(V)<\infty,\ \forall\,2\leq p<2+\frac{s-n}{d-s}.
\end{equation}
This is sharp for $s>d-1$ and has applications on Furstenberg sets and discretized sum-product problems.

Estimates \eqref{Lp-Orponen} and \eqref{DOV} are closely related, due to Orponen's formula 
\begin{equation}\label{Orponen-formula}
	\int\|\pi^y\mu\|_{L^p(G(d,n))}^p\,d\nu(y)=\int\|\pi_V\mu\|_{L^{p}(\pi_V\nu)}^{p}\,d\gamma_{d,n}(V)
\end{equation}
given $\mu$ has continuous density (see Lemma 3.1 in \cite{Orp19} and Lemma 4.17 in \cite{DOV22}), where
$$\pi^y\mu(V):=\int_{y+V^\perp}\mu\,d\mathcal{H}^{d-n}=\pi_V\mu(\pi_Vy).$$

Although the proof of \eqref{Orponen-formula} is not hard (just change variables $\R^d=V\oplus V^{\perp}$), it is very important. In \cite{Orp19}, the proof of \eqref{Lp-Orponen} starts from \eqref{Orponen-formula}; in \cite{DOV22}, \eqref{Orponen-formula} is used to obtain an incidence estimate from \eqref{DOV}.

In \cite{DOV22}, the authors suggest studying the mixed-norm
$$\int\|\pi_V\mu\|_{L^{p}(\mathcal{H}^n)}^{q}\,d\gamma_{d,n}(V),$$
especially for $q\leq p$. Their motivation is the following. As a corollary of \eqref{DOV}, for almost all $V\in G(d,n)$, 
\begin{equation}
	\label{range-p-DOV}
	\pi_V\mu\in L^p(\mathcal{H}^n|_V),\  \forall\, p<2+\frac{s-n}{d-s}.
\end{equation}
Since both $\supp\mu$, $G(d,n)$ are compact, one may expect a wider range of $p$ by considering smaller $q$.

In fact it is not hard to obtain mixed-norm estimates for $q\leq p$. It follows directly from Sobolev embedding that
\begin{equation}
	\label{Sobolev-embedding}
	\int\|\pi_V\mu\|_{L^p(\mathcal{H}^n)}^2\,d\gamma_{d,n}(V)<\infty, \forall\, p\begin{cases}
		<\frac{2n}{2n-s}, &n<s\leq 2n\\=\infty, &s> 2n
	\end{cases}.
\end{equation}
The range of $p$ in \eqref{Sobolev-embedding} is wider than \eqref{range-p-DOV} if $n<s<2(d-n)$. Then one can apply a mixed-norm version of the Riesz-Thorin interpolation to obtain 
$$\int\|\pi_V\mu\|_{L^{p}(\mathcal{H}^n)}^{q}\,d\gamma_{d,n}(V)<\infty$$
for $2\leq q\leq p<\infty$ and
\begin{equation}
	\label{interpolate-Sobolev-DOV}
	\frac{n}{p}+\frac{2d-2n-s}{q}>d-s.
\end{equation}

However, unfortunately, interpolation on $p$ never makes its range wider. In this paper we shall interpolate also on dimensions of measures, that indeed extends the range of $p$.

Now we turn to the mixed-norm
$$\int\|\pi^y\mu\|^q_{L^p(G(d,n))}d\nu(y).$$
When $q>p$ it is easily reduced to $q=p$ by considering
$$\int\|\pi^y\mu\|^p_{L^p(G(d,n))}\,f(y)d\nu(y),\ \|f\|_{L^{(q/p)'}(\nu)}=1;$$
while when $q<p$, due to the lack of Sobolev embedding, there seems no way to reduce it to $q=p$. As we mentioned above, the first step of Orponen's argument on \eqref{Lp-Orponen} is to apply his formula \eqref{Orponen-formula}, whose proof relies on $q=p$.
So we need new ideas, even to start.

To state our main theorem, we need the classical $s$-energy
$$I_s(\mu):=\iint|x-y|^{-s}\,d\mu(x)\,d\mu(y),$$
as well as a new quantity $A_s(\mu)$ called the $s$-amplitude. This new quantity plays a key role in our discussion, and makes our statement quite clean.

\begin{definition}\label{def-amplitude}
For every compactly supported Radon measure $\mu$ on $\R^d$ and $0<s<d$, define the $s$-amplitude of $\mu$ by
$$ A_s(\mu):=\left\||\cdot|^{-s}*\mu\right\|_{L^\infty}=\sup_{y\in\R^d}\int|x-y|^{-s}\,d\mu(x).$$	
\end{definition}

\begin{thm}
	\label{main-full-version}
Suppose $\mu, \nu$ are compactly supported Radon measures on $\R^d$, $0<t<n$, $0<s, \alpha<d$, and $s+t\geq 2n$. Denote
	$$q_0:=1+\frac{s+t-2n}{2(d-\alpha)}.$$
	Then, for every
	$$1\leq p< \frac{2nq_0}{n+t}=\frac{2n}{n+t}\cdot \left(1+\frac{s+t-2n}{2(d-\alpha)}\right),$$
	we have
	\begin{equation}
	\int\|\pi^y\mu\|^{q}_{L^{p}(G(d,n))}d\nu(y)\lesssim_{d,n,p,s,t,\alpha} \begin{cases}
	I_{s}(\mu)^{1/2}\cdot A_\alpha(\mu)^{q-1}\cdot I_t(\nu)^{1/2}, & q=q_0\\
	I_{s}(\mu)^{1/2}\cdot A_\alpha(\mu)^{q-1}\cdot A_{\max\{t, \frac{q}{2q_0}t\}}(\nu), & q>q_0
	\end{cases}.\end{equation}

	Furthermore, when $t=n$,
	\begin{equation}
		\int\|\pi^y\mu\|^{q}_{L^{q_0}(G(d,n))}d\nu(y)\lesssim_{d,n,s,\alpha}\begin{cases}
		I_{s}(\mu)^{1/2}\cdot A_\alpha(\mu)^{q-1}\cdot I_n(\nu)^{1/2}, & q=q_0\\
			I_{s}(\mu)^{1/2}\cdot A_\alpha(\mu)^{q-1}\cdot A_{\max\{n, \frac{q}{2q_0}n\}}(\nu), & q>q_0
		\end{cases}.
	\end{equation}
\end{thm}

Theorem \ref{main-full-version} implies the following mixed-norm estimates on Frostman measures.
\begin{corollary}
	\label{cor-p-q}
	Suppose $\mu, \nu$ are compactly supported Radon measures on $\R^d$, satisfying 
	$$\mu(B(x,r))\lesssim r^{s_\mu},\ \nu(B(x,r))\lesssim r^{s_\nu},\ \forall\,x\in\R^d,\,r>0,$$
	where $$0<s_\nu<n,\quad 2n-s_\nu<s_\mu<d.$$
Then 
	\begin{equation}
		\int\|\pi^y\mu\|^{q}_{L^{p}(G(d,n))}d\nu(y)<\infty
	\end{equation}
	if one of the following holds.
	\begin{itemize}
		\item $s_\mu\leq 2d-3n$ and 
$$1\leq p<\frac{2n}{n+\max\{2n-s_\mu, 0\}}\left(1+\frac{\max\{s_\mu-2n, 0\}}{2(d-s_\mu)}\right)=\begin{cases}
	\frac{2n}{3n-s_\mu},& s_\mu<2n\\2+\frac{s_\mu-2n}{d-s_\mu}, &s_\mu\geq 2n
\end{cases},$$ $$1\leq q< \frac{2s_\nu}{\max\{2n-s_\mu, 0\}}\left(1+\frac{\max\{s_\mu-2n, 0\}}{2(d-s_\mu)}\right)=\begin{cases}
	\frac{2s_\nu}{2n-s_\mu},& s_\mu<2n\\\infty, &s_\mu\geq 2n
\end{cases}.$$

\item $s_\mu\geq 2d-2n$ and 
$$1\leq p<\frac{2n}{n+s_\nu}\left(1+\frac{s_\mu+s_\nu-2n}{2(d-s_\mu)}\right),\quad 1\leq q< 2+\frac{s_\mu+s_\nu-2n}{d-s_\mu}.$$

\item $2d-3n<s_\mu< 2d-2n$ and $p,q$ lie in the region enclosed by 
$$1\leq p<\frac{2n}{n+s_\nu}\left(1+\frac{s_\mu+s_\nu-2n}{2(d-s_\mu)}\right),$$ $$1\leq q<\frac{2s_\nu}{\max\{2n-s_\mu, 0\}}\left(1+\frac{\max\{s_\mu-2n, 0\}}{2(d-s_\mu)}\right)=\begin{cases}
	\frac{2s_\nu}{2n-s_\mu},& s_\mu<2n\\\infty, &s_\mu\geq 2n
\end{cases},$$
and
$$\frac{s_\mu-2d+3n}{1-\frac{d-s_\mu}{n}p}< n+\frac{2d-2n-s_\mu}{\frac{d-s_\mu}{s_\nu}q-1}.$$
	\end{itemize}

In particular, in the case $p=q$,
$$\int\|\pi^y\mu\|^{p}_{L^{p}(G(d,n))}d\nu(y)<\infty$$
for every
$$1\leq p<\begin{cases}
		\frac{2n}{n+s_\nu}\cdot\left(1+\frac{s_\mu+s_\nu-2n}{2(d-s_\mu)}\right),&\text{if }s_\mu\geq 2d-3n\\
		\frac{2n}{n+\max\{2n-s_\mu, 0\}}\cdot\left(1+\frac{\max\{s_\mu-2n, 0\}}{2(d-s_\mu)}\right), &\text{if }s_\mu<2d-3n
	\end{cases}.$$
\end{corollary}
As a remark, in Corollary \ref{cor-p-q} $p,q\rightarrow\infty$ as $s_\mu\rightarrow d$, so it covers all $p,q\in[1,\infty)$. Also when $n=d-1$ and $p=q$, our result significantly improves \eqref{Lp-Orponen}. 

The assumption $0<t<n$ is necessary (see Section \ref{sec-alternative-DOV}), so Corollary \ref{cor-p-q} is not comparable to \eqref{DOV} by taking $\mu=\nu$ and $p=q$. But our method does give an alternative proof of \eqref{DOV}. See also Section \ref{sec-alternative-DOV}.

By ignoring $q$, we have the following range of $p$ for $\pi^y\mu(V)\in L^p(G(d,n))$.
\begin{corollary}
	\label{cor-p}
	Suppose $\mu, \nu$ are compactly supported Radon measures on $\R^d$, satisfying 
	$$\mu(B(x,r))\lesssim r^{s_\mu},\ \nu(B(x,r))\lesssim r^{s_\nu},\ \forall\,x\in\R^d,\,r>0,$$
	where $$0<s_\nu<n,\quad 2n-s_\nu<s_\mu<d.$$
	Then for every
	$$1\leq p<\begin{cases}
		\frac{2n}{n+s_\nu}\cdot\left(1+\frac{s_\mu+s_\nu-2n}{2(d-s_\mu)}\right),&\text{if }s_\mu\geq 2d-3n\\
		\frac{2n}{n+\max\{2n-s_\mu, 0\}}\cdot\left(1+\frac{\max\{s_\mu-2n, 0\}}{2(d-s_\mu)}\right), &\text{if }s_\mu<2d-3n
	\end{cases},$$
	there exists $y\in\supp\nu$ such that $$\pi^y\mu(V)\in L^p(G(d,n)).$$

	In particular, if $s_\mu>2n$, then set of $y\in\R^d$ such that
	$$\inf\{p: \pi^y\mu(V)\notin L^p(G(d,n))\}<
        2+\frac{s_\mu-2n}{d-s_\mu}$$
	has Hausdorff dimension $0$.
\end{corollary}

Notice that, near the critical line segment
$$\{(s_\mu, s_\nu)\in(0,d)^2: s_\mu+s_\nu=2n,\, 0<s_\nu<n\},$$
the critical $p$ in Corollary \ref{cor-p} equals $\frac{2n}{n+s_\nu}=\frac{2n}{3n-s_\mu}>1$. On the other hand, by embedding $\R^{n+1}$ into $\R^d$, the sharpness of the visibility problem in $\R^{n+1}$ implies that $s_\mu+s_\nu>2n$ is necessary for the existence of $p>1$. Together it follows that the range of $p$ has jump discontinuities at critical cases. See Figure \ref{fig:range-p} below: in the shadow $p>\frac{2n}{3n-s_\mu}$, while in the blank $p=1$. 
This is unexpected and so surprising. I think there are deep reasons behind this phenomenon to be explored.

\begin{figure}[H]
\centering
\begin{tikzpicture}[scale=0.75, p2/.style={line width=0.275mm, black}, p3/.style={line width=0.15mm, black!50!white}]

\fill[color=black!30] (3.5,0.5) -- (3.5,2) -- (2, 2);
\draw[->] (0, 0) -- (0, 3);
\draw[->] (0, 0) -- (5, 0);
\draw[-] (2, 2) -- (3.5, 0.5);
\draw[-] (3.5,0) -- (3.5, 2) -- (0,2);

\fill (2,2) circle (0.07);
\fill (3,1) circle (0.07);
\fill (3.5,0.5) circle (0.07);

\draw (2, 1.5) node[anchor=north east]{$p=1$};
\draw (2.75, 1) node[anchor=south west]{$p>\frac{2n}{3n-s_\mu}$};
\draw (5, 0) node[anchor=north]{$s_\mu$};
\draw (0, 3) node[anchor=east]{$s_\nu$};
\draw (3.5, 0) node[anchor=north]{$d$};
\draw (0, 2) node[anchor=east]{$n$};
\draw (3.5, 0.5) node[anchor=west]{$(d,2n-d)$};
\draw (2, 2) node[anchor=south]{$(n,n)$};
\draw (3.5, -1.5) node[anchor=east]{$n>d/2$};

\filldraw[fill=black!30] (13,0) -- (15,0) -- (15, 1.5) -- (11.5, 1.5);
\draw[->] (10, 0) -- (10, 3);
\draw[->] (10, 0) -- (16, 0);
\draw[-] (11.5, 1.5) -- (13, 0);
\draw[-] (15,0) -- (15, 1.5) -- (10,1.5);

\fill (11.5,1.5) circle (0.07);
\fill (12.5,0.5) circle (0.07);
\fill (13,0) circle (0.07);

\draw (10, 0.25) node[anchor=south west]{$p=1$};
\draw (12.25, 0.4) node[anchor=south west]{$p>\frac{2n}{3n-s_\mu}$};
\draw (16, 0) node[anchor=north]{$s_\mu$};
\draw (10, 3) node[anchor=east]{$s_\nu$};
\draw (13, 0) node[anchor=north]{$2n$};
\draw (10, 1.5) node[anchor=east]{$n$};
\draw (14, -1.5) node[anchor=east]{$n<d/2$};
\draw (15, 0) node[anchor=north]{$d$};
\draw (11.5, 1.5) node[anchor=south]{$(n,n)$};

\end{tikzpicture}
\caption{jump discontinuities of $p$ at $s_\mu+s_\nu=2n$}
\label{fig:range-p}
\end{figure}
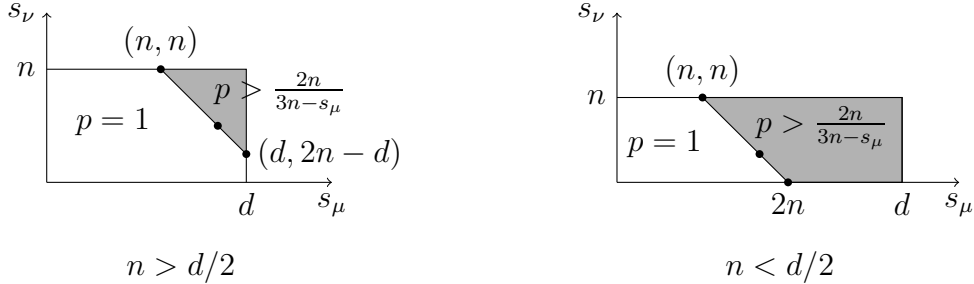

I do not know whether Theorem \ref{cor-p-q} and Theorem \ref{cor-p} are sharp, or how far they are from being sharp. It seems very hard to compute $L^p$-estimates from known geometric examples. I think the first step towards the sharpness is to understand the jump discontinuities along the critical line segment.

Our results also generalize the visibility problem. Notice that not every generalization is nontrivial. For example, one can easily conclude that for every $E, F\subset\R^d$, $\dH E>n$, $\dH E+\dH F>2n$, there exists $y\in F$ such that
\begin{equation}
	\label{trivial-generalization}
	\gamma_{d,d-n}\{W\in G(d,d-n):W\cap (E-y)\neq\emptyset\}>0.
\end{equation}
To see this, just project $E, F$ onto a $(n+1)$-dimensional subspace, preserving their dimensions, then \eqref{trivial-generalization} follows from the visibility problem in $\R^{n+1}$. See \cite{DIOWZ20} for an application of this trick on Falconer distance conjecture.

To make nontrivial generalizations, we consider the set of $m$-planes determined by $E-y$, that is,
$$\pi^y(E^m):=\{W\in G(d,m):W=\Span\{x_1-y,\dots,x_m-y\}: x_1,\dots,x_m\in E\}.$$
Throughout this paper, when writing
$$W=\Span\{\vec{v_1},\dots,\vec{v_m}\}$$
vectors $\vec{v_1},\dots,\vec{v_m}$ are a priori assumed nonzero and linearly independent.
\begin{corollary}\label{cor-geometric-consequence}
Suppose $1\leq m\leq d-1$ and $E\subset\R^d$.
\begin{enumerate}[(i)]
	\item If
	$$\dH E>\max\left\{d-\frac{m}{2m-1}, \ 3m-d\right\},$$
	then $$\dH\{y: \gamma_{d,m}(\pi^y(E^m))=0\}$$ $$\leq \max\left\{2(d-m)-\dH E,\ \frac{(d-m)(\dH E-2m+m(d-\dH E))}{d-m-m(d-\dH E)},\ 0\right\}.$$

	\item If
	$$\dH E>\max\left\{2(d-m), \ d-\frac{2m-d}{m-1}\right\},$$
	then $$\dH\{y: \gamma_{d,m}(\pi^y(E^m))=0\}=0.$$
\end{enumerate}
\end{corollary}

One can check that $(i)$ matches the results above on the visibility problem with $m=1$. Also $(ii)$ makes sense only if $m>d/2$.

To prove results like Corollary \ref{cor-geometric-consequence}, people used to consider multi-linear estimates. We shall see in Section \ref{sec-proof-thm-cor} that our interpolation between multi-linear estimates gives a better bound than multi-linear estimates themselves. We hope this will shed lights on other multiple configuration problems in the future.

As above I do not know how sharp Corollary \ref{cor-geometric-consequence} is or what to expect. One could follow the idea in \cite{Orp18} on the sharpness of the visibility problem, to transfer it to orthogonal projections via projective transformations. But then one should check, not only the projected image has positive Lebesgue measure, but also each slice is not contained in a lower dimensional affine subspace. Again we need examples.

{\bf Notation.} 

$X\lesssim Y$ means $X\leq CY$ for some constant $C>0$. $X\approx Y$ means $X\lesssim Y$ and $Y\lesssim X$. $X\lesssim_\epsilon Y$ means $X\leq CY$ for some constant $C=C(\epsilon)>0$.


\section{Preliminaries on Analytic Interpolation}
\subsection{}

Suppose $f$ is an analytic function on $\C$ such that
$$M_0:=\sup_{\Re z=0}|f(z)|<\infty,\quad M_1:=\sup_{\Re z=1}|f(z)|<\infty.$$
Hadamard three-lines lemma states that for every $z\in\{0<\Re z<1\}$,
$$|f(z)|\leq M_0^{\Re z}\cdot M_1^{1-\Re z}.$$

There is a higher dimensional version of Hadamard three-lines lemma. That is, given real numbers $a_1,\dots,a_k\in\R^n$ and an analytic function $f$ on $\C^n$, if $|f(z)|\leq 1$ for every $z\in a_j+i\R$, $\forall j$, then $|f(z)|\leq 1$ if $\Re z$ lies in the convex hull generated by $a_1,\dots,a_k$.

Hadamard three-lines lemma was first announced in 1890s. It was introduced into harmonic analysis in mid-20th century to bound norms of operators between $L^p$ spaces. These techniques are still widely used today, including Riesz-Thorin interpolation, Stein's interpolation, and many others.

\subsection{}

Our analytic interpolation is inspired by both classical ones and the work of Greenleaf-Iosevich \cite{GI12}, Grafakos-Greenleaf-Iosevich-Palsson \cite{GGIP15}. Suppose $\phi\in C_0^\infty$, $\phi\geq 0$, $\int \phi =1$, $\phi_\delta(\cdot):=\delta^{-d}\phi(\frac{\cdot}{\delta})$, and let
$$\mu^\delta(x):=\phi_\delta*\mu(x).$$
Then the Riesz potential,
\begin{equation}
	\label{Riesz-potential}
	\frac{\pi^\frac{z}{2}}{\Gamma(\frac{z}{2})}\,|\cdot|^{-d+z}*\mu^\delta(x),
\end{equation}
is a smooth function in $x\in\R^d$, initially defined for $\Re z>0$ and can be extended to $z\in\C$ by analytic continuation. It is well known that \eqref{Riesz-potential} has Fourier transform
$$\frac{\pi^\frac{d-z}{2}}{\Gamma(\frac{d-z}{2})}\,\hat{\mu^\delta}(\xi)|\xi|^{-z},$$
and in particular \eqref{Riesz-potential} equals $c_d\,\mu^\delta$ when $z=0$. We refer to \cite{GS77}, p71, p192, for details.

In \cite{GI12}\cite{GGIP15}, geometric configuration problems are reduced to
$$\Lambda(\mu^\delta,\cdots,\mu^\delta),$$
where $\mu$ is a Frostman measure and $\Lambda$ is a symmetric $k$-linear form. Take
\begin{equation}
	\label{Alex-analytic-multi-linear-form}
	\Phi(z_1,\dots,z_k):=\Lambda\left(\frac{\pi^{z_1/2}}{\Gamma(z_1/2)}\,|\cdot|^{-d+z_1}*\mu^\delta,\cdots,\frac{\pi^{z_k/2}}{\Gamma(z_k/2)}\,|\cdot|^{-d+z_k}*\mu^\delta\right)
\end{equation}
as an analytic function on $\C^k$. If $\mu(B(x,r))\lesssim r^{s+\epsilon}$ and $\Re z_j=d-s\in(0,d)$, it follows immediately that
\begin{equation}
	\label{positive-z}
	\left|\int |x-y|^{-d+z_j}\,d\mu(y)\right|\leq \int |x-y|^{-s}\,d\mu(y)\lesssim 1.
\end{equation}
Therefore $|\Phi(z)|$ can be reduced to a bilinear form 
$$B\left(\frac{\pi^{z_{j_1}/2}}{\Gamma(z_{j_1}/2)}\,|\cdot|^{-d+z_{j_1}}*\mu^\delta,\ \frac{\pi^{z_{j_2}/2}}{\Gamma(z_{j_2}/2)}\,|\cdot|^{-d+z_{j_2}}*\mu^\delta\right),$$
where
$$\Re z_{j_1}=\Re z_{j_2}=-\frac{k(d-s)}{2}.$$
Since $\Lambda$ is symmetric, the bilinear form estimate is independent in the choice of $j_1, j_2$, which means the same estimate holds for $\binom{k}{2}$ non-proportional vectors $(\Re z_1,\dots, \Re z_k)$ whose convex contains the origin. Hence estimates of
$$\Phi(0)= \Lambda(\mu^\delta,\cdots,\mu^\delta)$$ 
follow from Hadamard three-lines lemma.

\subsection{}

There are some technical issues about \cite{GI12}\cite{GGIP15} to be clarified.
\begin{enumerate}[(i)]
	\item Though \eqref{positive-z} holds, one cannot conclude that the Riesz potential \eqref{Riesz-potential} is also $\lesssim 1$, due to the unbounded factor $|\Gamma(z/2)|^{-1}$. To resolve this issue, we shall work with
\begin{equation}\label{def-mu-z}
	\mu^\delta_z(x):=e^{z^2}\frac{\pi^\frac{z}{2}}{\Gamma(\frac{z}{2})}\,|\cdot|^{-d+z}*\mu^\delta(x).
\end{equation}
The role of $e^{z^2}$ is to control
$$|e^{z^2}\cdot\Gamma^{-1}(z)|\lesssim_{\Re z} 1,$$
as $\Gamma^{-1}$ is an entire function of order $1$.

\medskip

\item When $\Re z\in(0,d)$, it is straightforward that
$$ |\mu^\delta_{z}(x)|\lesssim_{\Re z} \int |x-y|^{-d+\Re z}\,d\mu^\delta(y)\approx \mu^\delta_{\Re z}(x).$$
But we are not convinced that this relation can be extended to general $z\in\C$. In fact the right hand side is not guaranteed positive because it is defined via integration by parts. In our argument below we are super careful when taking absolute values of $\mu^\delta_z$, as we do not know how to take Fourier transform of $|\mu^\delta_z|$ to obtain
$$\hat{\mu^\delta}(\xi)|\xi|^{-\Re z}.$$
Readers can keep an eye on our timing of taking absolute values in Section \ref{sec-q=k}.
\end{enumerate}

\section{Energy, Amplitude, and Dimensions of Measures}\label{sec-def-amplitude}

For every $E\subset\R^d$, denote by $\mathcal{M}(E)$ the collection of compactly supported Radon measures on $E$.

The well-known Frostman Lemma implies that, for every Euclidean subset $E\subset\R^d$ and every $s<\dH E$, there exists $\mu\in\M(E)$ such that
$$\mu(B(x,r))\lesssim r^s,\ \forall\,r>0,\,\forall\,x\in\R^d.$$
In fact
$$\dH E=\sup\{s:\exists\, \mu\in\M(E):\sup_x \frac{\mu(B(x,r))}{r^s}<\infty \}.$$
We call $$\sup_x \frac{\mu(B(x,r))}{r^s}$$
the Frostman constant of $\mu$ (of dimension $s$).

By direct computation, the above implies that for every $s<\dH E$ there exists $\mu\in\M(E)$ such that the $s$-energy
$$I_s(\mu):=\iint|x-y|^{-s}\,d\mu(x)d\,\mu(y)$$
is finite. Also
\begin{equation}
	\label{dim-I-s}
	\dH E=\sup\{s:\exists\, \mu\in\M(E): I_s(\mu)<\infty\}.
\end{equation}

In dimension theories $s$-energy plays an important role due to its Fourier-analytic representation
$$I_s(\mu)=C_{d,s}\int|\hat{\mu}(\xi)|^2\,|\xi|^{-d+s}\,d\xi=C_{d,s}\,\|(-\Delta)^{-\frac{d-s}{4}}\mu\|_{L^2}^2.$$
In fact Kaufman's simple alternative proof of Marstrand projection theorem is just
\begin{equation}
	\begin{aligned}
		\int_{S^1} \|\pi_e\mu\|_{L^2(\R)}^2\,d\sigma(e)=&\iint|\widehat{\pi_e\mu}(r)|^2\,dr\,d\sigma(e)\\= &\iint|\hat{\mu}(re)|^2\,dr\,d\sigma(e)\\=&\int|\hat{\mu}(\xi)|^2|\xi|^{-1}\,d\xi\\= & C\,I_1(\mu).
	\end{aligned}
\end{equation}

So far, in the literature, all estimates on Frostman measures are written in terms of the Frostman constant and the energy . But these are not enough for our analytic interpolation. For our use, estimates should hold with complex-valued $\mu$, or more precisely, $\mu^\delta_z$ defined in \eqref{def-mu-z}. The $s$-energy works well with $\mu^\delta_z$: thanks to its Fourier-analytic representation,
\begin{equation}
	\label{I-s-mu-z}
	\begin{aligned}
		I_s(\mu^\delta_z)=&C_{d,s}\int|\widehat{\mu^\delta_z}(\xi)|^2\,|\xi|^{-d+s}\,d\xi\\ = &C_{d,s}\cdot \left|\frac{e^{z^2}\pi^\frac{d-z}{2}}{\Gamma(\frac{d-z}{2})}\right|^2\cdot \int|\hat{\mu}(\xi)|^2\,|\hat{\phi}(\delta\xi)|^2\,|\xi|^{-d+s-2\Re z}\,d\xi\\\leq &C_{d,s,\Re z}\cdot I_{s-2\Re z}(\mu).
	\end{aligned}
\end{equation}

However, there seems no easy way to deal with the Frostman constant of $\mu^\delta_z$, namely
$$\sup_{x}\frac{|\mu^\delta_z(B(x,r))|}{r^s}.$$
This is why we introduce the $s$-amplitude defined in Definition \ref{def-amplitude}. This definition is very natural. In fact, \eqref{dim-I-s}, the connection between $\dH E$ and $I_s(\mu)$, is built upon $A_s(\mu)<\infty$ (see, for example, \cite{Mat15}, Theorem 2.8), that immediately implies
$$\dH E=\sup\{s: \exists\,\mu\in\mathcal{M}(E),\  A_s(\mu)<\infty\}.$$
However there seems no further discussion about this quantity in the literature.

Unlike Frostman constant, $s$-amplitude can be easily extended to act on $\mu^\delta_z$: if $0<s-\Re z<d$, then
\begin{equation}
	\label{A-s-mu-z}
	\begin{aligned}
		A_s(\mu^\delta_z)=&\frac{\Gamma(\frac{d-s}{2})}{\pi^\frac{d-s}{2}}\cdot\left\|\frac{\pi^\frac{d-s}{2}}{\Gamma(\frac{d-s}{2})}|\cdot|^{-s}*\mu^\delta_z\right\|_{L^\infty}\\=& \left|\frac{e^{z^2}\pi^{\frac{d+z}{2}}\cdot\Gamma(\frac{d-s}{2})\cdot\Gamma(\frac{s-z}{2})}{\Gamma(\frac{s}{2})\cdot\Gamma(\frac{d-z}{2})\cdot\Gamma(\frac{d-(s-z)}{2})}\right|\cdot\left\||\cdot|^{-s+z}*\mu^\delta\right\|_{L^\infty}\\\leq & \,C_{d, s, \Re z}\cdot A_{s-\Re z}(\mu^\delta).
	\end{aligned}
\end{equation}
Notice that 
$$\left|\frac{\pi^\frac{d-s}{2}}{\Gamma(\frac{d-s}{2})}\right|\cdot A_s$$
can be extended to $s\in \C$. Therefore, if $0<\Re(w+z)<d$, then \eqref{A-s-mu-z} implies
\begin{equation}
	\label{sup-mu-z-w}
	\|(\mu^\delta_z)_w\|_{L^\infty}=\left|\frac{e^{w^2}\pi^\frac{w}{2}}{\Gamma(\frac{w}{2})}\right|\cdot A_{d-w}(\mu^\delta_z)\leq C_{d, \Re z, \Re w}\cdot A_{d-\Re (w+z)}(\mu^\delta),\ \forall\,z,w\in\C.
\end{equation}

It is routine to consider $\mu^\delta$ first and take $\delta\rightarrow 0$ at the very end. One can check that all implicit constants below are independent in $\delta$. From now we write $\mu$ for $\mu^\delta$ for abbreviation and assume $\mu$ has continuous density.

\section{Trivial Estimates on Orthogonal Projections}\label{sec-trivial}
\subsection{} 
The finiteness for $p=q=1$ is trivial. In fact for every finite measure $\mu$ we have
$$\int\pi^y\mu(V)\,d\gamma_{d,n}(V)\approx 1,\ \forall\,y\notin\supp\mu.$$
To see this, since $\gamma_{d,n}$ is induced by the haar measure $\theta_d$ on the orthogonal group $O(d)$, by its invariance
\begin{equation}\label{take-use-of-invariance-physical}
	\begin{aligned}
		\int \pi^y\mu(V)\,d\gamma_{d,n}(V)= &\int_{O(d)}\int_{\R^n}\mu(y+g\cdot(x',0))\,dx'\,d\theta_d(g)\\= &\int_{O(d)}\int_{S^{n-1}}\int_0^\infty\mu(y+g\cdot(r\sigma,0))\,r^{n-1}dr\,d\sigma\,d\theta_d(g)\\=&\int_{S^{n-1}}\int_0^\infty\left(\int_{O(d)}\mu(y+g\cdot(r\sigma,0))\,d\theta_d(g)\right)r^{n-1}dr\,d\sigma\\= &\frac{|S^{n-1}|}{|S^{d-1}|}\cdot\int_{S^{d-1}}\int\mu(y+r\sigma)\,r^{n-1}dr\,d\sigma\\= &\frac{|S^{n-1}|}{|S^{d-1}|}\cdot\iint\mu(y+x)|x|^{n-d}dx\approx 1.
	\end{aligned}
\end{equation}

As a consequence,
$$\iint\pi^y\mu(V)\,d\gamma_{d,n}(V)\,d\nu(y)\approx 1$$
if $\supp\mu$, $\supp\nu$ are disjoint. 

\subsection{}
The trivial estimate \eqref{take-use-of-invariance-physical} looks perfect and there seems nothing more to discuss. But here we would like to present another trivial estimate of 
$$\iint\pi^y\mu(V)\,d\gamma_{d,n}(V)\,d\nu(y)$$
that inspires our argument below.

Assume both $\mu, \nu$ have continuous density, then by Orponen's formula \eqref{Orponen-formula} we have, with $s+t=2n$,
\begin{equation}\label{both-positive}
\begin{aligned}
	\iint \pi^y\mu(V) \,d\gamma_{d,n}(V)\,d\nu(y)=&\iint\pi_V\mu(u)\cdot\pi_V\nu(u)\,d\mathcal{H}^{n}(u)\,d\gamma_{d,n}(V)\\\leq & \|(-\Delta)^{-\frac{n-s}{4}}\pi_V\mu\|_{L^2}\cdot \|(-\Delta)^{-\frac{n-t}{4}}\pi_V\nu\|_{L^2}\\= &C\cdot I^{1/2}_s(\mu)\cdot I^{1/2}_t(\nu).
\end{aligned}
\end{equation}
Here the second line of \eqref{both-positive} follows because in general
\begin{equation}
	\label{Plancherel}
	\int f\,g= \int \hat{f}\,\hat{g} =\int \hat{f}(\xi)|\xi|^{\alpha}\ \hat{g}(\xi)|\xi|^{-\alpha}\,d\xi\leq\|(-\Delta)^{\alpha/2}f\|_{L^2}\cdot\|(-\Delta)^{-\alpha/2}g\|_{L^2}.
\end{equation}
The last line of \eqref{both-positive} follows because for arbitrary but fixed $V_0\in G(d,n)$,
\begin{equation}\label{take-use-of-invariance-frequency}
\begin{aligned}
	\iint |\widehat{\pi_V\mu}(\xi)|^2\,|\xi|^{\alpha}\,d\xi\,d\gamma_{d,n}(V)=&\iint |\widehat{\pi_{gV_0}\mu}(\xi)|^2\,|\xi|^{\alpha}\,d\xi\,d\theta_d(g)\\=&\int\int_{gV_0} |\hat{\mu}(\xi)|^2\,|\xi|^{\alpha}\,d\mathcal{H}^n(\xi)\,d\theta_d(g)\\=&\int\int_{V_0} |\hat{\mu}(g\xi)|^2\,|\xi|^{\alpha}\,d\mathcal{H}^n(\xi)\,d\theta_d(g),
\end{aligned}
\end{equation}
which equals a constant multiple of $I_{n+\alpha}(\mu)$ by applying polar coordinate on $V_0$ and integrate in $\theta_d$ first as in \eqref{take-use-of-invariance-physical}.

The estimate \eqref{both-positive} looks useless, as we already know it is $\approx 1$ even if the energy blows up. However, it inspires our proof in Section \ref{sec-q=1}, which is the beginning of this project. We hope that presenting this ``trivial" estimate here would help readers have a better understanding of this paper.

\subsection{}
Before finally going to the proof, we would like to deal with a small technical issue. When $\mu$ is compactly supported of continuous density,
$$\pi_V\mu(u)=\int_{\pi_V^{-1}(u)}\mu\,d\mathcal{H}^{d-n}$$
and therefore
$$|\pi_V\mu(u)|\leq \mathcal{H}^{d-n}(\pi_V^{-1}(u)\cap\supp\mu)\cdot\|\mu\|_{L^\infty}\lesssim\|\mu\|_{L^\infty}.$$

However, in this paper we work with $\mu_z$, whose support is the whole space. To fix this we set
\begin{equation}\label{our-def-pi-V}
	\pi_V\mu(u):=\int_{\pi_V^{-1}(u)}\mu\,\psi\,d\mathcal{H}^{d-n},
\end{equation}
for some $\psi\in C_0^\infty(\R^d)$. Then
$$|\pi_V\mu_z(u)|\lesssim_\psi\|\mu_z\|_{L^\infty}.$$
Notice $\pi_V$ coincides with the classical definition when $\mu$ has compact support and $\psi=1$ on $\supp\mu$, so it does change anything in the main theorem. This modified $\pi_V$ brings a new issue, that is, now \eqref{take-use-of-invariance-frequency} ends up with $I_{n+\alpha}(\psi\mu_z)$, and it seems not as straightforward as \eqref{I-s-mu-z} to conclude
$$|I_s(\psi\mu_z)|\lesssim_{d,s,\Re z} I_{s-2\Re z}(\mu).$$
It is not hard, as one can easily see
$$|I_s(\psi\mu_z)|\lesssim_\psi I_s(\mu_z).$$
More precisely, by the Fourier-analytic representation of $s$-energy,
$$I_s(\psi\mu_z)=\int \left|\int\hat{\mu_z}(\eta)\,\hat{\psi}(\xi-\eta)\,d\eta\right|^2|\xi|^{-d+s}\,d\xi\lesssim_\psi \iint |\hat{\mu_z}(\eta)|^2\,|\hat{\psi}(\xi-\eta)|\,d\eta\,|\xi|^{-d+s}\,d\xi.$$
When $|\xi|>|\eta|/2$, 
$$\int_{|\xi|>|\eta|/2} |\hat{\psi}(\xi-\eta)|\,|\xi|^{-d+s}\,d\xi\lesssim |\eta|^{-d+s}\int |\hat{\psi}(\xi-\eta)|\,d\xi\lesssim |\eta|^{-d+s}.$$
When $|\xi|<|\eta|/2$, we have $|\xi-\eta|\gtrsim |\eta|$ and therefore
$$\int_{|\xi|<|\eta|/2}|\hat{\psi}(\xi-\eta)|\,|\xi|^{-d+s}\,d\xi\lesssim_N\int_{|\xi-\eta|\gtrsim |\eta|}|\xi-\eta|^{-N}\,|\xi|^{-d+s}\,d\xi\lesssim_N (1+|\eta|)^{s-N},$$
where the last inequality follows from changing variables $\xi=|\eta|\zeta$.

From now $\pi_V\mu$ is defined as \eqref{our-def-pi-V}.

\section{$q=1$}\label{sec-q=1}
We start our proof with $q=1$. This is also how we discover jumps of $p$ at $s+t=2n$.
\begin{prop}\label{prop-q=1}
Suppose $\mu$ is a complex-valued smooth function on $\R^d$, $\nu$ is a compactly supported measure on $\R^d$, and $0<t<n$. Then
	\begin{equation}
		\label{estimate-q=1}
		\int\|\pi^y\mu\|_{L^p(G(d,n))}d\nu(y)\lesssim_{d,n,p,t} I_{2n-t}(\mu)^{1/2}\cdot I_{t}(\nu)^{1/2}
	\end{equation}
	for every
	$$1\leq p<\frac{2n}{n+t}.$$
	Furthermore,
	\begin{equation}
		\label{estimate-q=1-t=n}
		\int\|\pi^y\mu\|_{L^1(G(d,n))}d\nu(y)\lesssim_{d,n} I_n(\mu)^{1/2}\cdot I_{n}(\nu)^{1/2}.
	\end{equation}
\end{prop}

By considering $\nu*\phi_\delta\rightarrow\nu$, we may assume $\nu$ has continuous density. The reason we work with complex-valued $\mu$ is for the analytic interpolation in Section \ref{sec-analytic-interpolation}. It brings extra difficulties and even the case $p=1$ is no longer trivial. More precisely, if we follow the argument in Section \ref{sec-trivial},
\begin{equation}\label{obstacle-complex-valued}
	\begin{aligned}
		\int\|\pi^y\mu\|_{L^1(G(d,n))}d\nu(y)=&\int\int_{G(d,n)} |\pi^y\mu(V)| \,d\gamma_{d,n}(V)\,d\nu(y)\\=&\iint|\pi_V\mu(u)|\cdot\pi_V\nu(u)\,d\mathcal{H}^{n}(u)\,d\gamma_{d,n}(V),
	\end{aligned}
\end{equation}
getting stuck due to the absolute value symbol. The following argument overcomes this obstacle, and eventually proves Proposition \ref{prop-q=1}.

Now we start our proof. The endpoint case $t=n$, $p=1$ follows directly from taking Cauchy-Schwarz of \eqref{obstacle-complex-valued}. So we assume $0<t<n$ and $1<p<\frac{2n}{n+t}<2$. 

To estimate $\|\pi^y\mu\|_{L^p(G(d,n))}$, $1\leq p<\infty$, the most natural way is to consider
$$\sup_{\|f\|_{L^{p'}}=1}\left|\int \pi^y\mu(V)\cdot f(V)\,d\gamma_{d,n}(V)\right|.$$
But this setup is not convenient to us due to the integral in $d\nu(y)$ outside supreme. Instead, we take $f$ to be the maximizer directly, namely
$$f_y(V)=f(V,y)=\frac{\sgn (\pi^y\mu(V))\cdot|\pi^y\mu(V)|^{p-1}}{\|\pi^y\mu\|_{L^p(G(d,n))}^{p-1}}$$
with
$$\sgn (\pi^y\mu(V)):=\frac{\pi^y\mu(V)}{|\pi^y\mu(V)|}\cdot \chi_{\pi^y\mu(V)\neq 0}(V,y).$$
Then the mixed-norm in Proposition \ref{prop-q=1} is reduced to
\begin{equation}\label{dual-form}
	\iint \pi^y\mu(V)\cdot f(V,y)\,d\gamma_{d,n}(V)\,d\nu(y)
\end{equation}
with $\|f(\cdot,y)\|_{L^{p'}}=1$, $\forall\,y$. 

We first fix $V$, parametrize $\R^d$ by $(u,v)\in V\oplus V^\perp$ and integrate in $v\in V^\perp$ first. Since $\pi^y\mu(V)=\pi_V\mu(\pi_Vy)$ is a constant in each level set $y\in\pi_V^{-1}(u)$, the integral \eqref{dual-form} equals
\begin{equation}
	\label{int-V-perp-first}
	\int_{G(d,n)}\int_{V} \pi_V\mu(u)\left(\int_{y\in \pi_V^{-1}(u)}\,f(V, y)\,\nu(y)\,d\mathcal{H}^{d-n}(y)\right)\,d\mathcal{H}^n(u)\,d\gamma_{d,n}(V).
\end{equation}

Denote
\begin{equation}
	\label{def-F}
	F(V, u):=\int_{y\in \pi_V^{-1}(u)}\,f(V, y)\,\nu(y)\,d\mathcal{H}^{d-n}(y).
\end{equation}
Then, similar to \eqref{both-positive}, the integral \eqref{int-V-perp-first} is reduced to
\begin{equation}\label{Cauchy-Schwarz}
		\int_{G(d,n)}\int_{V} \pi_V\mu(u)\,F(V,u)\,d\mathcal{H}^n(u)\,d\gamma_{d,n}(V)\leq I^{1/2}_{2n-t}(\mu)\cdot \|(-\Delta)^{-\frac{n-t}{4}}F\|_{L^2(\mathcal{H}^n\times\gamma_{d,n})}.
\end{equation}

The first factor is desirable, so it remains to consider the second. Since $0<t<n$, it follows that
\begin{equation}\label{back-to-physical-q=1}
	\begin{aligned}
		\|(-\Delta)^{-\frac{n-t}{4}}F\|_{L^2}^2= & \int\int_{\R^n} |\hat{F}(V,\xi)|^2\,|\xi|^{-n+t}\,d\xi\,d\gamma_{d,n}(V)\\=&C_{d,n,s}\iiint |u-u'|^{-t}\,F(V, u)\,F(V, u')\,du\,du'\,d\gamma_{d,n}(V).
	\end{aligned}
\end{equation}
By the definition of $F$ in \eqref{def-F} and again $\R^d=V\oplus V^{\perp}$, \eqref{back-to-physical-q=1} equals, up to a multiplicative constant,
\begin{equation}
	\label{back-to-f-g}
	\iiint |\pi_V(y-y')|^{-t}\,f(V, y)\,d\nu(y)\,f(V, y')\,d\nu(y')\,d\gamma_{d,n}(V).
\end{equation}

Now fix $y\neq y'$ and integrate over $G(d,n)$ first. Since $1<p<2$ and $\|f(\cdot,y)\|_{L^{p'}}=1$ for every $y$, it follows that $$\|f(\cdot, y)f(\cdot, y')\|_{L^{p'/2}}=1,\ \forall\,y, y',$$ with $1<p'/2<\infty$. Therefore by H\"older's inequality 
\begin{equation}
	\label{int-V}
	\int|\pi_V(y-y')|^{-t}\,f(V, y)\,f(V, y')\,d\gamma_{d,n}(V)\leq \left(\int|\pi_V(y-y')|^{-t\cdot (p'/2)'}\,d\gamma_{d,n}(V)\right)^{\frac{1}{(p'/2)'}}.
\end{equation}
It is well known (see, for example, Theorem 3.12 in \cite{Mat95}) that \eqref{int-V} is integrable and 
$$\lesssim_{d,n,p,s} |y-y'|^{-t}$$
if $$t\cdot (p'/2)'<n.$$ 
After computation it is equivalent to
$$p<\frac{2n}{n+t}$$
which completes the proof of Proposition of \ref{prop-q=1}. 

\section{$q=k$}\label{sec-q=k}
Proposition \ref{prop-q=1} only works when $p<2$. This is not satisfactory. In this section we prove results for large $p$. For technical reasons both $\mu, \nu$ have to be positive. See \eqref{L-infty-mu-z} in the proof below.

\begin{prop}
	\label{prop-q=k}
	Suppose $\mu, \nu$ are compactly supported Radon measures on $\R^d$, $k\geq 2$, $0<t<n$, and $0<s, \alpha<d$ are real numbers satisfying
    \begin{equation}\label{k-s-t}
    s+t=2n+2(k-1)(d-\alpha).
    \end{equation}
	Then
	\begin{equation}\label{estimate-q=k}
		\int\|\pi^y\mu\|^{k}_{L^p(G(d,n))}d\nu(y)\lesssim_{d,n,p, k, s,t} I_s(\mu)^{1/2}\cdot A_\alpha(\mu)^{k-1}\cdot I_t(\nu)^{1/2}
	\end{equation}
	for every
	$$1\leq p<\frac{2nk}{n+t}.$$

	Furthermore, when $t=n$,
	\begin{equation}
		\label{estimate-q=k-t=n}
		\int\|\pi^y\mu\|^{k}_{L^{k}(G(d,n))}d\nu(y)\lesssim_{d,n,k,s} I_{s}(\mu)^{1/2}\cdot A_\alpha(\mu)^{k-1}\cdot I_{n}(\nu)^{1/2}.
	\end{equation}
\end{prop}

As $G(d,n)$ is compact, we many assume $p\geq k$ and denote $p_k:=p/k\geq 1$. Write
$$\int\|\pi^y\mu\|^k_{L^p(G(d,n))}d\nu(y)=\int\|(\pi^y\mu)^k\|_{L^{p_k}(G(d,n))}d\nu(y).$$
 We shall study the $p_k$-norm of the $k$-linear form $\prod_{j=1}^k\pi^y\mu_j$, with $\mu_j=\mu$. 

If one repeats the argument in Section \ref{sec-q=1}, there is an obstacle that \eqref{Plancherel} does not work for multi-linear forms. To overcome this difficulty, we take ideas from Greenleaf-Iosevich \cite{GI12} and Grafakos-Greenleaf-Iosevich-Palsson \cite{GGIP15}, that are already sketched in Section \ref{sec-def-amplitude}. Compared with their multi-linear estimates, our case is more subtle because of the mixed-norm. We should be super careful, especially on the timing of taking absolute values.

Now we start the proof of Proposition \ref{prop-q=k}. Again we assume both $\mu, \nu$ have continuous density. With $\mu_j=\mu$, $j=1,\dots,k$, similar to Section \ref{sec-q=1} it suffices to consider
\begin{equation}\label{dual-k-form}
	\iint \prod_{j=1}^k\pi^y\mu_j(V)\cdot f(V,y)\,d\gamma_{d,n}(V)\,d\nu(y),
\end{equation}
with $\|f(\cdot,y)\|_{L^{p_k'}}=1$, $\forall\,y$. 

Let
$$\Phi(z_1,\dots,z_k):=\iint \prod_{j=1}^k\pi^y\mu_{z_j}(V)\cdot f(V,y)\,d\gamma_{d,n}(V)\,d\nu(y)$$
be an analytic function on $\C^k$, where
$$\mu_z(x):=\frac{e^{z^2}\pi^\frac{z}{2}}{\Gamma(\frac{z}{2})}\,|\cdot|^{-d+z}*\mu(x)$$
as in \eqref{def-mu-z}. Similar to Section \ref{sec-q=1}, since $\pi_V\mu_{z}$ is a constant in each level set $\pi_V^{-1}(u)$, one can write $\Phi(z_1,\dots,z_k)$ as
$$\int_{G(d,n)}\int_{V} \pi_V\mu_{z_{j_0}}(u)\left(\int_{y\in \pi_V^{-1}(u)}\prod_{j\neq j_0}\pi^y\mu_{z_j}(V)\,f(V, y)\,\nu(y)\,d\mathcal{H}^{d-n}(y)\right)\,d\mathcal{H}^n(u)\,d\gamma_{d,n}(V),$$
for arbitrary $j_0\in\{1,\dots,k\}$. Take $$\Re z_j=d-\alpha,\ j\neq j_0, \text{ and }\Re z_{j_0}=-\sum_{j\neq j_0}\Re z_j=-(k-1)(d-\alpha).$$ 
Since $\mu$ is positive,
\begin{equation}\label{L-infty-mu-z}
	|\mu_{z_j}(x)|=\left|\frac{e^{z_j^2}\pi^\frac{z_j}{2}}{\Gamma(\frac{z_j}{2})}\int|x-y|^{-d+z_j}d\mu(y)\right|\lesssim_{d, \alpha} A_\alpha(\mu), \ \forall\,j\neq j_0.
\end{equation}

Denote
\begin{equation}
	\label{def-F-k}
	F(V, u):=\int_{y\in \pi_V^{-1}(u)}\prod_{j\neq j_0}\pi^y\mu_{z_j}(V)\,f(V, y)\,\nu(y)\,d\mathcal{H}^{d-n}(y).
\end{equation}

Then, as in \eqref{both-positive} and \eqref{Cauchy-Schwarz},
\begin{equation}
	\label{Cauchy-Schwarz-q=k}
	|\Phi(z_1,\dots,z_k)|=\left| \iint \pi_V\mu_{z_{j_0}}(u)\,F(V,u)\,du\,dV\right|\leq I_{2n-t}(\mu_{z_{j_0}})^{1/2}\cdot \|(-\Delta)^{-\frac{n-t}{4}}F\|_{L^2}.
\end{equation}

By \eqref{k-s-t}, we have
$$2n-t=s-2(k-1)(d-\alpha)=s+2\Re z_{j_0},$$
so by \eqref{I-s-mu-z} the square of the first factor is
$$I_{2n-t}(\mu_{z_{j_0}})\lesssim I_{2n-t-2\Re z_{j_0}}(\mu)=I_{s}(\mu),$$
as desired. 

It remains to estimate the second factor in \eqref{Cauchy-Schwarz-q=k}. Now we can feel free to take absolute values. By \eqref{def-F-k}, \eqref{L-infty-mu-z},
\begin{equation}
	\label{abs-F}
	|F(V, u)|\lesssim_{d,\alpha} A_\alpha(\mu)^{k-1}\cdot \int_{y\in \pi_V^{-1}(u)}|f(V, y)|\,\nu(y)\,d\mathcal{H}^{d-n}(y).
\end{equation}

Then the special case \eqref{estimate-q=k-t=n} follows quickly: when $t=n$, $p_k=1$, we have $\|f(\cdot,y)\|_{L^\infty}=1$, $\forall\,y$, and therefore
$$\|F\|_{L^2(\mathcal{H}^n\times\gamma_{d,n})}\lesssim A_\alpha(\mu)^{k-1}\cdot\|\pi_V\nu\|_{L^2(\mathcal{H}^n\times\gamma_{d,n})}\lesssim A_\alpha(\mu)^{k-1}\cdot I_{n}(\nu)^{1/2},$$
as desired.

We proceed to prove \eqref{estimate-q=k}. When $0<t<n$, as in Section \ref{sec-q=1} we write $\|(-\Delta)^{-\frac{n-t}{4}}F\|_{L^2}^2$ to be
\begin{equation}
	\label{back-to-physical-q=k}
	\int\int_{\R^n} |\hat{F}(V,\xi)|^2\,|\xi|^{-n+t}\,d\xi\,d\gamma_{d,n}(V)=\iiint |u-u'|^{-t}\,F(V, u)\overline{F(V, u')}\,du\,du'\,dV.
\end{equation}

By plugging \eqref{abs-F} into \eqref{back-to-physical-q=k}, we obtain
$$A_\alpha^{2(k-1)}(\mu)\cdot \iiint |\pi_V(y-y')|^{-t}\,|f(V, y)|\,d\nu(y)\,|f(V, y')|\,d\nu(y')\,d\gamma_{d,n}(V).$$

The rest is the same as Section \ref{sec-q=1}, and eventually gets to
$$p/k=p_k<\frac{2n}{n+t}.$$

As $j_0\in\{1,\dots,k\}$ is arbitrary, the estimate holds for $k$ non-proportional vectors $(\Re z_1,\dots, \Re z_k)\in\R^k$ whose sum equals $\vec{0}$. It implies the origin lies in their convex hull. Hence by Hadamard three-lines lemma
$$ |\Phi(0)|=\left|\iint \prod_{j=1}^k\pi^y\mu_j(V)\cdot f(V,y)\,d\gamma_{d,n}(V)\,d\nu(y)\right|\lesssim I_s(\mu)^{1/2}\cdot A_\alpha(\mu)^{k-1}\cdot I_t(\nu)^{1/2},$$
that completes the proof of Proposition \ref{prop-q=k}.

\section{Analytic Interpolation}\label{sec-analytic-interpolation}
Previous sections already constitute an interesting paper. One can even prove a weaker version of Corollary \ref{cor-geometric-consequence} from Proposition \ref{prop-q=k} as a geometric application. But, if we compute the range of $p$ for
$$\pi^y\mu(V)\in L^p(G(d,n)),$$
there is something strange. Suppose $\mu, \nu$ are Frostman measures satisfying
$$\mu(B(x,r))\lesssim r^{s_\mu},\ \nu(B(x,r))\lesssim r^{s_\nu},\ \forall\,x\in\R^d, \forall\,r>0,$$
with $0<s_\nu<n$ and $2n-s_\nu<s_\mu<2n$. It follows from Proposition \ref{prop-q=1}, \ref{prop-q=k} that
$$\pi^y\mu(V)\in L^p(G(d,n)), \text{ for some }y\in\supp\nu,$$
for every
\begin{equation}
	\label{p-max-k}
	1\leq p<\max_k\frac{2nk}{3n-s_\mu+2(k-1)(d-s_\mu)},
\end{equation}
where $k$ is taken over all positive integers satisfying
\begin{equation}
	\label{s-mu-s-nu-k}
	s_\mu+t_\nu\geq 2n+2(k-1)(d-s_\mu).
\end{equation}
By solving $k$ from \eqref{p-max-k}, we see that when $s_\mu\leq 2d-3n$ we should take $k$ as small as possible, that is $k=1$, while when $s_\mu>2d-3n$ we should take $k$ as large as possible, that is, by \eqref{s-mu-s-nu-k},
$$k=1+\left[\frac{s_\mu+s_\nu-2n}{2(d-s_\mu)}\right],$$
where $[\cdot]$ denotes the integer part. This means the range of $p$ has jump discontinuities at each $\frac{s_\mu+s_\nu-2n}{2(d-s_\mu)}\in\Z_+$. We have seen in the introduction that jumps at $s_\mu+s_\nu-2n=0$ do exist. However, there are no reasons to support jumps elsewhere. 

We would like to drop the $[\cdot]$ symbol for a wider range of $p$. As we comment right after \eqref{interpolate-Sobolev-DOV}, traditional interpolations only on $p,q$ do not help. To make it, we introduce a new technique that also interpolates dimensions of measures. For technical reasons, one cannot interpolate between estimates in Proposition \ref{prop-q=k} directly. We have pointed out that both $\mu, \nu$ there have to be positive. In fact for the proof of Theorem \ref{main-full-version} readers can skip Section \ref{sec-q=k} and read this section directly. However, without Proposition \ref{prop-q=k} and the observation above, there is no way to know what to prove. For later use of this mechanism, one can first obtain an analog of Proposition \ref{prop-q=k} on the scratch paper, and then prove desired estimates by arguments below.

We shall see how $s$-amplitude helps us, and why we consider $s, \alpha$ to be possibly different in Proposition \ref{prop-q=k}.

\begin{prop}
	\label{prop-p-q}
	Suppose $\mu, \nu$ are compactly supported Radon measures on $\R^d$, $0<t<n$, $0<s, \alpha<d$ and $s+t\geq 2n$. Denote
	\begin{equation}
		\label{relation-q-s-t-alpha}
		q:=1+\frac{s+t-2n}{2(d-\alpha)}.
	\end{equation}
	Then
	\begin{equation}
		\label{estimate-p-q}
		\int\|\pi^y\mu\|^{q}_{L^{p}(G(d,n))}d\nu(y)\lesssim_{d,n,p,s,t,\alpha} I_{s}^{1/2}(\mu)\cdot A_\alpha(\mu)^{q-1}\cdot I_t(\nu)^{1/2}
	\end{equation}
	for every
	\begin{equation}
		\label{range-p-general}
		1\leq p< \frac{2nq}{n+t}=\frac{2n}{n+t}\cdot \left(1+\frac{s+t-2n}{2(d-\alpha)}\right).
	\end{equation}

	Furthermore, when $t=n$,
	\begin{equation}
		\label{estimate-p-q-t=n}
		\int\|\pi^y\mu\|^{q}_{L^{q}(G(d,n))}d\nu(y)\lesssim_{d,n,s,\alpha} I_{s}(\mu)^{1/2}\cdot A_\alpha(\mu)^{q-1}\cdot I_{n}(\nu)^{1/2}.
	\end{equation}
\end{prop}

The case $q=1$ is already done in Proposition \ref{prop-q=1}, so we assume $p,q>1$.

It suffices to consider
	$$\iint \pi^y\mu(V)\cdot f(V,y)\,d\gamma_{d,n}(V)\,g(y)\,d\nu(y),$$
with $\|f(\cdot,y)\|_{L^{p'}}=1$, $\forall\,y$, and $\|g\|_{L^{q'}(\nu)}=1$. 

For every $z\in\C$, let
$$s_z=s+2z,$$
$$\alpha_z=\alpha+z,$$
$$q_z=1+\frac{s_z+t-2n}{2(d-\alpha_z)}=\frac{2(d-\alpha)+s+t-2n}{2(d-\alpha-z)},$$
$$p_z=\frac{p}{q}\cdot q_z.$$
Then we make the following observations:
\begin{itemize}
	\item $s_0=s$, $\alpha_0=\alpha$, $q_0=q$, $p_0=p$;
	\item for every $z\in\C$, \begin{equation}
	\label{k-s-t-z}
	s_z+t=2n+2(q_z-1)(d-\alpha_z);
          \end{equation}
    \item both $\frac{1}{p_z}, \frac{1}{q_z}$ are linear in $z$.
\end{itemize}

Now take
$$\Phi(z)=\iint \pi^y\mu_z(V)\cdot f_z(V,y)\,d\gamma_{d,n}(V)\,g_z(y)\,d\nu(y)$$
as an analytic function, where $\mu_z$ is defined as \eqref{def-mu-z},
$$f_z(V,y):=\sgn(f)\cdot |f(V,y)|^{p'\left(1-\frac{1}{p_z}\right)},$$
$$g_z(y):=\sgn(g)\cdot|g(y)|^{q'\cdot\left(1-\frac{1}{q_z}\right)}.$$

Since both $\frac{1}{p_z}, \frac{1}{q_z}$ are linear in $z$,
\begin{equation}
	\label{f-z-L-p}
	\begin{aligned}
		\|f_z(\cdot, y)\|_{L^{p_{\Re z}'}}=&\left\||f(\cdot, y)|^{p'\left(1-\frac{1}{p_{\Re z}}\right)}\right\|_{L^{p_{\Re z}'}}=1,\ \forall\,y;\\\|g_z\|_{L^{q'_{\Re z}}}=&\left\||g|^{q'\left(1-\frac{1}{q_{\Re z}}\right)}\right\|_{L^{q'_{\Re z}}}=1.
	\end{aligned}
\end{equation}

Heuristically, if Proposition \ref{prop-q=k} could be applied to $\mu_z$, it would imply that for each $q_{\Re z}\in\Z_+$,
\begin{equation}
	\label{invoke-estimate-at-integer}
	\begin{aligned}
		|\Phi(z)|^{q_{\Re z}}\lesssim &\,I_{\Re s_z}^{1/2}(\mu_z)\cdot A_{\Re \alpha_z}(\mu_z)^{q_{\Re z}-1}\cdot I_t(\nu)^{1/2}\\\lesssim &\,I_{\Re (s_z-2z)}(\mu)^{1/2}\cdot A_{\Re(\alpha_z-z)}(\mu)^{q_{\Re z}-1}\cdot I_t(\nu)^{1/2}\\=& \,I_{s}(\mu)^{1/2}\cdot A_\alpha(\mu)^{q_{\Re z}-1}\cdot I_t(\nu)^{1/2}.
	\end{aligned}
\end{equation}
Then the desired estimate of $|\Phi(0)|$ would follow from Hadamard three-lines lemma.

In the following we make the estimate \eqref{invoke-estimate-at-integer} rigorous. 

First consider $\Re z=-\frac{s+t-2n}{2}<0$. In this case $q_{\Re z}=1$. Then by H\"older's inequality, \eqref{f-z-L-p}, Proposition \ref{prop-q=1}, and \eqref{I-s-mu-z}, 
$$|\Phi(z)|\leq \int \|\pi^y\mu_z\|_{L^{p_{\Re z}}(G(d,n))}d\nu(y)\lesssim I_{\Re s_z}(\mu_z)^{1/2}\cdot I_t(\nu)^{1/2}\lesssim_{d, s, \Re z} I_s(\mu)^{1/2}\cdot I_t(\nu)^{1/2}.$$

We also need estimates with $\Re z>0$. Since $q_{\Re z}\rightarrow\infty$ as $\Re z\rightarrow d-\alpha>0$, there exists $\Re z>0$ such that $q_{\Re z}$ is an integer $k\geq 2$. In this case, by H\"older's inequality and \eqref{f-z-L-p},
$$|\Phi(z)|^k\leq \int \|\pi^y\mu_z\|^k_{L^{p_{\Re z}}(G(d,n))}d\nu(y)=\iint \pi^y\mu_z(V)\cdots \pi^y\mu_z(V)\cdot h(V,y)\,d\gamma_{d,n}(V)\,d\nu(y),$$
where $$\|h(\cdot,y)\|_{L^{\left(\frac{p_{\Re z}}{k}\right)'}}=1, \ \forall\,y. $$

Fix $z$ and let
$$\Psi_z(w_1,\dots,w_k):=\iint \prod_{j=1}^k\pi^y(\mu_z)_{w_j}(V)\cdot h(V,y)\,d\gamma_{d,n}(V)\,d\nu(y)$$
be an analytic function on $\C^k$. For an arbitrary but fixed $1\leq j_0\leq k$, take
$$\Re w_j=d-\Re \alpha_z,\ j\neq j_0, \text{ and }\Re w_{j_0}=-\sum_{j\neq j_0}\Re z_j=-(k-1)(d-\Re \alpha_z).$$
Then, with
$$F(V, u):=\int_{y\in \pi_V^{-1}(u)}\prod_{j\neq j_0}\pi^y(\mu_z)_{w_j}(V)\,h(V, y)\,\nu(y)\,d\mathcal{H}^{d-n}(y),$$
$|\Psi_z(w_1,\dots,w_k)|$ can be written as
$$\begin{aligned}
	&\left|\int_{G(d,n)}\int_{V} \pi_V(\mu_z)_{w_{j_0}}(u)\cdot F(V,u)\,d\mathcal{H}^n(u)\,d\gamma_{d,n}(V)\right|\\\leq \ &I_{2n-t}((\mu_z)_{w_{j_0}})^{1/2}\cdot \|(-\Delta)^{-\frac{n-t}{4}}F\|_{L^2}\\\lesssim \ &I_{2n-t-2\Re(z+w_{j_0})}(\mu)^{1/2}\cdot \|(-\Delta)^{-\frac{n-t}{4}}F\|_{L^2}.
\end{aligned}
$$

Due to $q_{\Re z}=k$, \eqref{k-s-t-z} and our choice of $w_{j_0}$, the first factor equals $I_s(\mu)^{1/2}$, as desired. 

To estimate the second factor, this time we need \eqref{sup-mu-z-w}. 

When $t=n$, it follows immediately that
$$\|F\|_{L^2(\mathcal{H}^n\times\gamma_{d,n})}\lesssim A_{\Re(\alpha_z-z)}(\mu)^{k-1}\cdot\|\pi_V\nu\|_{L^2(\mathcal{H}^n\times\gamma_{d,n})}\lesssim A_\alpha(\mu)^{k-1}\cdot I_{n}(\nu)^{1/2},$$
as desired. 

When $0<t<n$,
$$\begin{aligned}
	&\ \|(-\Delta)^{-\frac{n-t}{4}}F\|^2_{L^2}\\=&
	\iiint |u-u'|^{-t}\,F(V, u)\overline{F(V, u')}\,du\,du'\,dV\\
	\lesssim &	\ A_{\Re (\alpha_z-z)}^{2(k-1)}(\mu)\cdot \iiint |\pi_V(y-y')|^{-t}\,|h(V, y)|\,d\nu(y)\,|h(V, y')|\,d\nu(y')\,d\gamma_{d,n}(V)\\=&\ A_{\alpha}^{2(k-1)}(\mu)\cdot \iiint |\pi_V(y-y')|^{-t}\,|h(V, y)|\,d\nu(y)\,|h(V, y')|\,d\nu(y')\,d\gamma_{d,n}(V).
\end{aligned}$$

The rest is the same as Section \ref{sec-q=1}. We omit details. Eventually it ends up with 
$$|\Psi_z(w_1,\dots,w_k)|\lesssim I_{s}(\mu)^{1/2}\cdot A_\alpha(\mu)^{k-1}\cdot I_t(\nu)^{1/2},$$
where the implicit constant is independent in the choice of $j_0$. By \eqref{I-s-mu-z} and \eqref{sup-mu-z-w}, this implicit constant is also independent in $\Im w$, $\Im z$. Therefore by Hadamard three-lines lemma, when $q_{\Re z}=k\geq 2$,
$$|\Phi(z)|^{q_{\Re z}}=|\Psi_z(0)|\lesssim I_{s}(\mu)^{1/2}\cdot A_\alpha(\mu)^{q_{\Re z}-1}\cdot I_t(\nu)^{1/2},$$
where the implicit constant is independent in $\Im z$.

Now we have estimates of $|\Phi(z)|$ at hand, for both positive and negative $\Re z$. Hence 
$$|\Phi(0)|^q\lesssim I_{s}(\mu)^{1/2}\cdot A_\alpha(\mu)^{q-1}\cdot I_t(\nu)^{1/2}$$
by Hadamard three-lines lemma again, which completes the proof of Proposition \ref{prop-p-q}.

\section{Proof of Corollary \ref{cor-p}, \ref{cor-geometric-consequence}: Optimize the Range of $p$}\label{sec-proof-thm-cor}

To finally go from Proposition \ref{prop-p-q} to Theorem \ref{main-full-version}, it remains to broaden the range of $q$. We leave it to the next section. In fact Proposition \ref{prop-p-q} is enough for Corollary \ref{cor-p}, \ref{cor-geometric-consequence}, as we should ignore $q$ for the maximal possible $p$.

\subsection*{Proof of Corollary \ref{cor-p}}

Since we need
$$q_0=1+\frac{s+t-2n}{2(d-\alpha)}\geq 1,$$
the assumption $s+t\geq 2n$ is required. Also for the finiteness of $I_s(\mu), A_\alpha(\mu), I_t(\nu)$, we need $0<s, \alpha<s_\mu$ and $0<t<s_\nu$. Therefore the maximal $p$ equals, up to the end point,
\begin{equation}
	\label{sup-s-alpha-t}
	\sup_{\substack{
	0\leq s,\alpha\leq s_\mu\\ 0\leq t\leq s_\nu\\ s+t\geq 2n}} \frac{2n}{n+t}\cdot \left(1+\frac{s+t-2n}{2(d-\alpha)}\right).
\end{equation}

The assumptions $0<s_\nu<n$ and $2n-s_\nu<s_\mu<d$ ensure the supreme is well defined, namely not taken over an empty set. As parameters lie in a compact set, the supreme can be attained, say at $s_0, \alpha_0, t_0$. 

Notice that, for every fixed $t$, the critical $p$ is an increasing function in $s, \alpha$. Therefore $s_0=\alpha_0=s_\mu$. To find $t_0$, write
$$\frac{2n}{n+t}\cdot \left(1+\frac{s_\mu+t-2n}{2(d-s_\mu)}\right)=\frac{n}{d-s_\mu}\left(1+\frac{2d-3n-s_\mu}{n+t}\right).$$
To make this quantity large, one can see that, when $s_\mu\geq 2d-3n$ we should take $t_0$ as large as possible, namely $t_0=s_\nu$; when $s_\mu<2d-3n$ we should take $t_0$ as small as possible, namely $t_0=\max\{0, 2n-s_\mu\}$. This completes the proof of Corollary \ref{cor-p}.

\subsection*{Proof of Corollary \ref{cor-geometric-consequence}}

When $m=1$ it matches previous result on the visibility problem, so assume $m\geq 2$.

Let $F\subset\R^d$, $\supp E\cap\supp F\neq\emptyset$, and $\mu, \nu$ be Frostman measures on $E, F$ of exponents $s_\mu, s_\nu$, respectively. It suffices to show that when $s_\nu$ is large enough then there exists $y\in F$ such that $\gamma_{d,m}(\pi^y(E))>0$. 

As in either case $\dH E>d-1$, there exist compact subsets $E_1,\dots, E_m\subset E$ such that $\mu(E_i)>0$, $\forall i$, and no $m$-tuple $(x_1,\dots, x_m)\in E_1\times\cdots\times E_m$ lies in a $m$-dimensional affine subspace. This means the map 
$$\Phi_y(x_1,\dots, x_m)=\Span\{x_1-y,\dots,x_m-y\}\in G(d, m)$$
is well defined on $E_1\times\cdots\times E_m$, for all $y\in F$. We claim that $\Phi_y$ has no critical point on $E_1\times\cdots\times E_m$. To see this, just consider the action of the orthogonal group on a neighborhood of $E_1\times\cdots\times E_m$. From this point of view, if there exists a critical point in this neighborhood, then every point in this neighborhood is critical, contradiction.

Denote $\mu_i=\mu|_{E_i}$ and define a measure $\Phi_y(\mu_1\times\cdots\times\mu_m)$ on $\pi^y(E^m)$ by
$$\int f(W)\,d\Phi_y(\mu_1\times\cdots\times\mu_m)(W)=\int\cdots\int f(\Phi_y(x_1,\dots, x_m))\,d\mu_1(x_1)\dots d\mu_m(x_m).$$
Then, to show $\gamma_{d,m}(\pi^y(E^m))>0$, it suffices to show the measure $\Phi_y(\mu_1\times\cdots\times\mu_m)$ has $L^p$ density for some $p>1$.

By considering $\mu*\phi_\delta\rightarrow\mu$, we may assume $\mu$ has continuous density. Since $\Phi_y$ is regular and $E_1\times\cdots\times E_m$ is compact, by the co-area formula we have
$$\Phi_y(\mu_1\times\cdots\times\mu_m)(W)\approx\int_{\Phi_y^{-1}(W)}\mu_1\cdots \mu_m\,d\mathcal{H}^{m^2}=\prod_{i=1}^m\int_{y+W}\mu_i\,d\mathcal{H}^m\leq |\pi^y\mu(W^\perp)|^m.$$

Therefore, to show $\Phi_y(\mu_1\times\cdots\times\mu_m)$ has $L^p$ density, it suffices to show
$$\int|\pi^y\mu(W^\perp)|^{mp}\,d\gamma_{d,m}(W)=\int|\pi^y\mu(V)|^{mp}\,d\gamma_{d,d-m}(V)<\infty.$$

To proceed, one could apply multi-linear estimates from Proposition \ref{prop-q=k}. However, as we already observed at the beginning of Section \ref{sec-analytic-interpolation}, the range of $p$ would have jumps. This means one has to make $s_\mu, s_\nu$ large enough to ensure the range of $p$ jumps across $m$. This explains why our interpolation between multi-linear estimates gives a better bound than multi-linear estimates themselves. We believe it will shed lights on other multiple configuration problems in the future.

To finish the proof we invoke Corollary \ref{cor-p} with $n=d-m$ and check when the range of $p$ covers $m=d-n$. Recall $0<s_\nu<n$ and $s_\mu+s_\nu>2n$ are always required.

We first prove $(i)$ in Corollary \ref{cor-geometric-consequence}. 

Since $s_\mu>2d-3n$, we need to slove
$$\frac{2n}{n+s_\nu}\cdot\left(1+\frac{s_\mu+s_\nu-2n}{2(d-s_\mu)}\right)>d-n.$$
This can be reduced to
\begin{equation}
	\label{solve-for-s-nu}
	(n+s_\nu)(n-(d-n)(d-s_\mu))>n(s_\mu-2d+3n).
\end{equation}
Since $n+s_\nu<2n$ is required, \eqref{solve-for-s-nu} has a solution only if 
$$n-(d-n)(d-s_\mu)>\frac{1}{2}(s_\mu-2d+3n)>0,$$ 
equivalent to
$$s_\mu>d-\frac{d-n}{2(d-n)-1}=d-\frac{m}{2m-1}.$$
Then we solve for $s_\nu$ from \eqref{solve-for-s-nu} to obtain
\begin{equation}
	\label{lower-bound-s-nu}
	s_\nu>\frac{n(s_\mu-2(d-n)+(d-n)(d-s_\mu))}{n-(d-n)(d-s_\mu)}=\frac{(d-m)(s_\mu-2m+m(d-s_\mu))}{d-m-m(d-s_\mu)},
\end{equation}
which completes the proof of $(i)$ in Corollary \ref{cor-geometric-consequence}. One can check that the right hand side of \eqref{lower-bound-s-nu} is non-negative unless
$$s_\mu>d-\frac{d-2n}{d-n-1}=d-\frac{2m-d}{m-1}.$$

Now we turn to $(ii)$ in Corollary \ref{cor-geometric-consequence}.

When $s_\mu>2d-3n$, it follows directly from taking the right hand side of $(i)$ to be $0$. So we assume $s_\mu<2d-3n$. Since $s_\mu>2n$, we need to solve
$$2+\frac{s_\mu-2n}{d-s_\mu}>d-n,$$
which is equivalent to 
$$s_\mu>d-\frac{d-2n}{d-n-1}=d-\frac{2m-d}{m-1},$$
that completes the proof.

There is one case we did not discuss, that is $s_\mu<2d-3n$ and $s_\mu<2n$. This is because in this case we need
$$\frac{2n}{3n-s_\mu}>d-n,$$
that has a solution only if $n=d-1$, namely $m=1$, already ruled out.

\section{Proof of Theorem \ref{main-full-version}, Corollary \ref{cor-p-q}: Broaden the Range of $q$}\label{sec-large-q}

\subsection*{Proof of Theorem \ref{main-full-version}}

Since both $G(d,n)$ and $\supp\nu$ are compact, for every pair $p, q$ in Proposition \ref{prop-p-q}, the estimate also holds for smaller $p,q$. So it remains to prove results for large $q$. Theorem \ref{main-full-version} is a summary of Proposition \ref{prop-p-q} and the following.
\begin{prop}
	\label{prop-large-q}
	Suppose $\mu, \nu$ are compactly supported Radon measures on $\R^d$, $0<t<n$, $0<\alpha<d$, and $2n-t\leq s<d$. Denote
	$$q_0:=1+\frac{s+t-2n}{2(d-\alpha)}.$$
	Then
	\begin{equation}\label{estimate-p-large-q}\int\|\pi^y\mu\|^{q}_{L^{p}(G(d,n))}d\nu(y)\lesssim_{d,n,p,s,t} I_{s}^{1/2}(\mu)\cdot A_\alpha(\mu)^{q-1}\cdot A_{\max\{t, \frac{q}{2q_0}t\}}(\nu)\end{equation}
	for every $q>q_0$ and
	$$1\leq p< \frac{2nq_0}{n+t}=\frac{2n}{n+t}\cdot \left(1+\frac{s+t-2n}{2(d-\alpha)}\right).$$

	Furthermore, when $t=n$,
	\begin{equation}
		\label{estimate-p-large-q-t=n}
		\int\|\pi^y\mu\|^{q}_{L^{q}(G(d,n))}d\nu(y)\lesssim_{d,n,q,s} I_{s}(\mu)^{1/2}\cdot A_\alpha(\mu)^{q-1}\cdot A_{\max\{t, \frac{q}{2q_0}t\}}(\nu).
	\end{equation}
\end{prop}

The proof goes by analyzing those right hand sides in Proposition \ref{prop-p-q}. We do not need to treat $t=n$ and $0<t<n$ separately. 

It suffices to consider
$$\int\|\pi^y\mu\|^{q_0}_{L^{p}(G(d,n))}\,g(y)d\nu(y),\ \|g\|_{L^{(q/q_0)'}(\nu)}=1.$$
When $q<2q_0$, we treat $g$ as $g\in L^2$. So we may assume $q\geq 2q_0$. By Proposition \ref{prop-p-q} it is bounded above by
$$I_{s}^{1/2}(\mu)\cdot A_\alpha(\mu)^{q-1}\cdot I_t(g\,d\nu)^{1/2}.$$
Now it remains to estimate
$$I_t(g\,d\nu)=\iint |y-y'|^{-t}\,g(y)\,d\nu(y)\,g(y')\,d\nu(y').$$
Let $r:=\frac{q}{2q_0}\geq 1$. Then the relation
$$\frac{1}{r}=1-\left(\frac{1}{(q/q_0)'}-\frac{1}{q/q_0}\right)$$
is satisfied. Denote
$$K(y,y')=|y-y'|^{-t}$$
as the kernel. By definition of $A_s$ we have
$$\left(\int |K(y,y')|^r\,d\nu(y)\right)^{1/r}=\left(\int |K(y,y')|^r\,d\nu(y)\right)^{1/r}\leq A_{\frac{q}{2q_0}t}(\nu)^{1/r}.$$
Then by H\"older's inequality, the fact $\|g\|_{L^{(q/q_0)'}(\nu)}=1$, and Young's inequality,
$$\iint K(y,y')\,g(y)\,d\nu(y)\,g(y')\,d\nu(y')\leq \left\|\int K(y,\cdot)\,g(y)\,d\nu(y)\right\|_{L^{q/q_0}(\nu)}\leq A_{\frac{q}{2q_0}t}(\nu)^{1/r},$$
as desired.

\subsection*{Proof of Corollary \ref{cor-p-q}}

It follows directly from Theorem \ref{main-full-version}. 

We first consider the special case $p=q$. Observe from Theorem \ref{main-full-version} that, lifting $q$ from $q_0$ to $2q_0$ does not change the finiteness. Then, since $$p<\frac{2nq_0}{n+t}\leq 2q_0,$$it follows that lifting $q$ from $q_0$ to $p$ does not influence the optimal $p$, that is Corollary \ref{cor-p}.

Now we consider general $p,q$. We need to find all $p,q\in[1,\infty)$ such that, there exist $s, \alpha, t$ satisfying 
$$0<s, \alpha<s_\mu,\ 0<\max\{t, \frac{q}{2q_0}t\}<s_\nu,\ s+t\geq 2n,\ p<\frac{2nq_0}{n+t},$$
where
$$q_0:=1+\frac{s+t-2n}{2(d-\alpha)}.$$

Since both $p,q$ are increasing in $s,\alpha$, up to the end point we may take $s=\alpha=s_\mu$ for large $p,q$. Then conditions above are equivalent to
$$\max\{2n-s_\mu, 0\}< t<s_\nu,$$
\begin{equation}
	\label{relation-p-q-t-s-mu}
	\begin{aligned}
p<\frac{2n}{n+t}\left(1+\frac{s_\mu+t-2n}{2(d-s_\mu)}\right)&=\frac{n}{d-s_\mu}\left(1+\frac{2d-3n-s_\mu}{n+t}\right),\\
q<\frac{2s_\nu}{t}\left(1+\frac{s_\mu+t-2n}{2(d-s_\mu)}\right)&=\frac{s_\nu}{d-s_\mu}\left(1+\frac{2d-2n-s_\mu}{t}\right).
	\end{aligned}
\end{equation}

When $s_\mu\leq 2d-3n$, both $p,q$ are decreasing in $t$, so for large $p,q$ we should take $t$ as small as possible, namely $t\searrow\max\{2n-s_\mu, 0\}$. Therefore the range of $p,q$ is
$$p<\frac{2n}{n+\max\{2n-s_\mu, 0\}}\left(1+\frac{\max\{s_\mu-2n, 0\}}{2(d-s_\mu)}\right)=\begin{cases}
	\frac{2n}{3n-s_\mu},& s_\mu<2n\\2+\frac{s_\mu-2n}{d-s_\mu}, &s_\mu\geq 2n
\end{cases},$$ $$q< \frac{2s_\nu}{\max\{2n-s_\mu, 0\}}\left(1+\frac{\max\{s_\mu-2n, 0\}}{2(d-s_\mu)}\right)=\begin{cases}
	\frac{2s_\nu}{2n-s_\mu},& s_\mu<2n\\\infty, &s_\mu\geq 2n
\end{cases}.$$

When $s_\mu\geq2d-2n$, both $p,q$ are increasing in $t$, so for large $p,q$ we should take $t$ as large as possible, that is $t\nearrow s_\nu$. Therefore the range of $p,q$ is
$$p<\frac{2n}{n+s_\nu}\left(1+\frac{s_\mu+s_\nu-2n}{2(d-s_\mu)}\right),\quad q< 2+\frac{s_\mu+s_\nu-2n}{d-s_\mu}.$$

When $2d-3n<s_\mu< 2d-2n$, $p$ is increasing in $t$ while $q$ is decreasing in $t$, so there is a relation between $p,q$. By solving for $t$ from \eqref{relation-p-q-t-s-mu} we end up with the region enclosed by
$$ 1\leq p<\frac{2n}{n+s_\nu}\left(1+\frac{s_\mu+s_\nu-2n}{2(d-s_\mu)}\right),$$ $$1\leq q<\frac{2s_\nu}{\max\{2n-s_\mu, 0\}}\left(1+\frac{\max\{s_\mu-2n, 0\}}{2(d-s_\mu)}\right)=\begin{cases}
	\frac{2s_\nu}{2n-s_\mu},& s_\mu<2n\\\infty, &s_\mu\geq 2n
\end{cases},$$
and
$$\frac{s_\mu-2d+3n}{1-\frac{d-s_\mu}{n}p}<n+\frac{2d-2n-s_\mu}{\frac{d-s_\mu}{s_\nu}q-1}.$$

\section{An Alternative Proof of Dabrowski-Orponen-Villa}\label{sec-alternative-DOV}
When $\mu=\nu$, $p=q$, Orponen's formula says
$$\int\|\pi^y\mu\|_{L^p(G(d,n))}^p\,d\nu(y)=\int\|\pi_V\mu\|_{L^{p+1}(\mathcal{H}^n)}^{p+1}\,d\gamma_{d,n}(V).$$
Then it is natural to ask whether Proposition \ref{prop-p-q} covers \eqref{DOV}. The answer is unfortunately no, due to the constraint $0<t<n$. In general this condition $0<t<n$ cannot be relaxed. Technically we need \eqref{back-to-physical-q=1}, \eqref{back-to-physical-q=k} to come back to the physical space and integrate in $V$. Geometrically, if our results hold for some $t>n$, then $s$ is allowed to be $<n$, contradicting the visibility problem in $\R^{n+1}$.

Despite this, it does not mean our method is not strong enough. In fact, when $\mu=\nu$ there are extra symmetries. With ideas from previous sections, it is straightforward to recover \eqref{DOV}.

\begin{prop}
	\label{prop-improve-DOV}
	Suppose $\mu$ is a compactly supported measure on $\R^d$, $0<s, \alpha<d$, and $2\leq q<\infty$ satisfying
		\begin{equation}
		\label{s-alpha}
		s=n+(q-2)(d-\alpha).
	\end{equation} 
	Then
	\begin{equation}
		\label{pi-V-q-q}
		\int\|\pi_V\mu\|_{L^{q}(\mathcal{H}^n)}^{q}\,d\gamma_{d,n}(V)\lesssim_{d,n,p} I_{s}(\mu)\cdot A_\alpha(\mu)^{q-2}.
	\end{equation}

	Moreover, for every $2\leq q< p\leq\infty$,
	\begin{equation}
		\label{pi-V-p-q}
		\int\|\pi_V\mu\|_{L^{p}(\mathcal{H}^n)}^{q}\,d\gamma_{d,n}(V)\lesssim_{d,n,p} I_{s}(\mu)\cdot A_\alpha(\mu)^{q-2},
	\end{equation}
	with
	\begin{equation}
		\label{s-alpha-p-q}
		s=n+(q-2)(d-\alpha)+n(1-\frac{q}{p}).
	\end{equation}
\end{prop}

When $s=\alpha$, the condition \eqref{s-alpha} becomes
$$q=2+\frac{s-n}{d-s}=\frac{2d-n-s}{d-s},$$
that coincides with \eqref{DOV}.

There are two ways to obtain \eqref{pi-V-p-q} from \eqref{pi-V-q-q}: run interpolation on $p,q$ with \eqref{Sobolev-embedding} as in the introduction, or invoke Sobolev embedding directly
\begin{equation}
	\begin{aligned}
		\int\|\pi_V\mu\|_{L^{p}(\mathcal{H}^n)}^{q}\,d\gamma_{d,n}(V)\lesssim &\int\|\pi_V(-\Delta)^{\frac{n}{2}(\frac{1}{q}-\frac{1}{p})}\mu\|_{L^{q}(\mathcal{H}^n)}^{q}\,d\gamma_{d,n}(V)
	\end{aligned}
\end{equation}
and then run the argument below. We leave details to readers. From now we only consider $p=q$.

As the idea is already explained clearly in Section \ref{sec-analytic-interpolation}, we decide to skip computations on the scratch paper and present the rigorous proof directly. Readers can follow the explanation in Section \ref{sec-analytic-interpolation} to figure out where the exponents below come from.

Now we start the proof. It suffices to consider
	$$\iint \pi_V\mu(u)\cdot f(V, u)\,d\mathcal{H}^n(u)\,d\gamma_{d,n}(V),\ \|f\|_{L^{q'}}=1.$$

For every $z\in\C$, let
$$s_z=s+2z,$$
$$\alpha_z=\alpha+z,$$
$$q_z=2+\frac{s_z-n}{d-\alpha_z}=\frac{2(d-\alpha)+s-n}{d-\alpha-z}.$$
Observe that
\begin{itemize}
	\item $s_0=s$, $\alpha_0=\alpha$, $q_0=q$;
	\item for every $z\in\C$, 
	\begin{equation}
	\label{s-alpha-z}
	s_z=n+(q_z-2)(d-\alpha_z);
\end{equation}
    \item $\frac{1}{q_z}$ is linear in $z$.
\end{itemize}

Now take
$$\Phi(z)=\iint \pi_V\mu_z(u)\cdot f_z(V,u)\,d\mathcal{H}^n(u)\,d\gamma_{d,n}(V)$$
as an analytic function, where $\mu_z$ is defined as \eqref{def-mu-z} and
$$f_z(V,y):=\sgn(f) \cdot |f(V,y)|^{q'\left(1-\frac{1}{q_z}\right)}.$$

Since $\frac{1}{q_z}$ is linear in $z$,
$$\|f_z\|_{L^{q_{\Re z}'}}=\left\||f|^{q'\left(1-\frac{1}{q_{\Re z}}\right)}\right\|_{L^{q_{\Re z}'}}=1.$$

When $\Re z=-\frac{s-n}{2}<0$, we have $q_{\Re z}=2$. Then by Cauchy-Schwarz, the classical $L^2$-estimate of orthogonal projections, \eqref{I-s-mu-z} and \eqref{s-alpha-z}, it follows that
$$|\Phi(z)|^2\leq \|\pi_V\mu_z\|^2_{L^2}=I_n(\mu_z)\lesssim I_{s}(\mu).$$

For positive $\Re z$, notice that $q_{\Re z}\rightarrow\infty$ as $\Re z\rightarrow d-\alpha>0$. Therefore there exists $\Re z>0$ such that $q_{\Re z}=2k\geq 4$ is an even integer. Then
$$|\Phi(z)|^{2k}\leq \int\|\pi_V\mu_z\|_{L^{2k}}^{2k}\,d\gamma_{d.n}(V)=\iint\pi_V\mu_z\cdot\pi_V\overline{\mu_z}\cdots\pi_V\mu_z\cdot\pi_V\overline{\mu_z}.$$

Fix this $z$ and let
$$\Psi_z(w_1,\dots,w_{2k}):= \iint\pi_V(\mu_z)_{w_1}\cdot\pi_V(\overline{\mu_z})_{w_2}\cdots\pi_V(\mu_z)_{w_{2k-1}}\cdot\pi_V(\overline{\mu_z})_{w_{2k}}$$
be an analytic function on $\C^{2k}$. For arbitrary but fixed $j_1, j_2\in\{1,\dots,2k\}$, take
$$\Re w_{j_1}=\Re w_{j_2}=-(k-1)(d-\Re \alpha_z), \ \text{and }\Re w_j=d-\Re \alpha_z,\ j\neq j_1, j_2.$$
Then, by \eqref{sup-mu-z-w}, for every $j\neq j_1, j_2$,
$$\|(\mu_z)_{w_j}\|_{L^\infty}\lesssim_{\Re z, \Re w_j} A_{d-\Re(z+w_j)}(\mu)=A_\alpha(\mu).$$
Therefore $|\Psi_z(w_1,\dots,w_{2k})|$ is bounded above by
$$\|\pi_V(\mu_z)_{w_{j_1}}\|_{L^2}^2\cdot A_{\alpha}^{2(k-1)}(\mu) \lesssim I_{n-2\Re(z+w_{j_1})}(\mu)\cdot A_\alpha^{2(k-1)}(\mu)=I_{s}(\mu)\cdot A_\alpha^{2(k-1)}(\mu).$$

As $j_1, j_2$ are arbitrary, by Hadamard three-lines lemma,
\begin{equation}\label{estimate-q=2k}
	|\Phi(z)|^{q_{\Re z}}=|\Psi_z(0)|\lesssim I_{s}(\mu)\cdot A_\alpha^{q_{\Re z}-2}(\mu),
\end{equation}
for every even integer $q_{\Re z}\geq 4$, where the implicit constant is independent in $\Im z$.

Now we have estimates of $|\Phi(z)|$ for both positive and negative $\Re z$, with implicit constants independent in $\Im z$. Hence by Hadamard three-lines lemma again,
$$|\Phi(0)|^{q}=|\Psi_z(0)|\lesssim I_{s}(\mu)\cdot A_\alpha^{q-2}(\mu),$$
that completes the proof.

\bibliographystyle{abbrv}
\bibliography{/Users/MacPro/Dropbox/Academic/paper/mybibtex.bib}

\begin{thebibliography}{10}

\bibitem{DOV22}
D.~Dabrowski, T.~Orponen, and M.~Villa.
\newblock Integrability of orthogonal projections, and applications to
  furstenberg sets.
\newblock {\em arXiv preprint arXiv:2107.04471}, 2022.

\bibitem{DIOWZ20}
X.~Du, A.~Iosevich, Y.~Ou, H.~Wang, and R.~Zhang.
\newblock An improved result for {F}alconer's distance set problem in even
  dimensions.
\newblock {\em Math. Ann.}, 380(3-4):1215--1231, 2021.

\bibitem{GS77}
I.~M. Gel'fand and G.~E. Shilov.
\newblock {\em Generalized functions. {V}ol. 1}.
\newblock Academic Press [Harcourt Brace Jovanovich, Publishers], New
  York-London, 1964 [1977].
\newblock Properties and operations, Translated from the Russian by Eugene
  Saletan.

\bibitem{GGIP15}
L.~Grafakos, A.~Greenleaf, A.~Iosevich, and E.~Palsson.
\newblock Multilinear generalized {R}adon transforms and point configurations.
\newblock {\em Forum Math.}, 27(4):2323--2360, 2015.

\bibitem{GI12}
A.~Greenleaf and A.~Iosevich.
\newblock On triangles determined by subsets of the {E}uclidean plane, the
  associated bilinear operators and applications to discrete geometry.
\newblock {\em Anal. PDE}, 5(2):397--409, 2012.

\bibitem{GIOW18}
L.~Guth, A.~Iosevich, Y.~Ou, and H.~Wang.
\newblock On {F}alconer's distance set problem in the plane.
\newblock {\em Invent. Math.}, 219(3):779--830, 2020.

\bibitem{Kau68}
R.~Kaufman.
\newblock On {H}ausdorff dimension of projections.
\newblock {\em Mathematika}, 15:153--155, 1968.

\bibitem{KS18}
T.~Keleti and P.~Shmerkin.
\newblock New {B}ounds on the {D}imensions of {P}lanar {D}istance {S}ets.
\newblock {\em Geom. Funct. Anal.}, 29(6):1886--1948, 2019.

\bibitem{Liu18}
B.~Liu.
\newblock An {$L^2$}-identity and pinned distance problem.
\newblock {\em Geom. Funct. Anal.}, 29(1):283--294, 2019.

\bibitem{Mar54}
J.~M. Marstrand.
\newblock Some fundamental geometrical properties of plane sets of fractional
  dimensions.
\newblock {\em Proc. London Math. Soc. (3)}, 4:257--302, 1954.

\bibitem{Mat87}
P.~Mattila.
\newblock Spherical averages of {F}ourier transforms of measures with finite
  energy; dimension of intersections and distance sets.
\newblock {\em Mathematika}, 34(2):207--228, 1987.

\bibitem{Mat95}
P.~Mattila.
\newblock {\em Geometry of Sets and Measures in Euclidean Spaces: Fractals and
  Rectifiability}, volume~44 of {\em Cambridge Studies in Advanced
  Mathematics}.
\newblock Cambridge University Press, Cambridge, 1995.

\bibitem{Mat15}
P.~Mattila.
\newblock {\em Fourier analysis and Hausdorff dimension}, volume 150.
\newblock Cambridge University Press, 2015.

\bibitem{Orp18}
T.~Orponen.
\newblock A sharp exceptional set estimate for visibility.
\newblock {\em Bull. Lond. Math. Soc.}, 50(1):1--6, 2018.

\bibitem{Orp19}
T.~Orponen.
\newblock On the dimension and smoothness of radial projections.
\newblock {\em Anal. PDE}, 12(5):1273--1294, 2019.

\bibitem{PS00}
Y.~Peres and W.~Schlag.
\newblock Smoothness of projections, {B}ernoulli convolutions, and the
  dimension of exceptions.
\newblock {\em Duke Math. J.}, 102(2):193--251, 2000.

\end{thebibliography}

\end{document}